\documentclass[11pt]{article}

\usepackage{amsfonts, amsbsy, amsmath, amssymb}

\usepackage{graphicx,float,hyper ref,listings}

\renewcommand{\to}{\rightarrow}

\newcommand{\Z}{\mathbb{Z}}

\newcommand{\F}{\mathbb{F}}

 \newcommand{\Sp}{\mathbb{S}}

\newcommand{\qed}{\hbox{\rule{6pt}{6pt}}}
\newcommand{\kc}[2]{ {\rm Col}_{#1}{(#2)}}
\newcommand{\nkc}[2]{ {\rm Col}^N_{#1}{(#2)}} 
\newcommand{\skc}[2]{ {\rm SCol}_{#1}{(#2)}}

\newcommand{\col}{ {\rm Col} }
\newcommand{\Lq}{ {\rm Lq} }
\newcommand{\mlq}{ {\rm MLq} }
\newcommand{\NI}{ {\rm NI} }
\newcommand{\bg}{ {\rm Bg} }

\newtheorem{theorem}{Theorem}[section]
\newtheorem{corollary}[theorem]{Corollary}
\newtheorem{lemma}[theorem]{Lemma}
\newtheorem{proposition}[theorem]{Proposition}

\newtheorem{definition}[theorem]{Definition}

\newtheorem{example}[theorem]{Example}

\newtheorem{remark}[theorem]{Remark}

\newtheorem{conjecture}[theorem]{Conjecture}

\newtheorem{appentable}[theorem]{Table}

\input{epsf.sty}
\usepackage{graphicx}

\begin{document}

\title{ Quandle Colorings of Knots  
and Applications}

\author{W. Edwin Clark, Mohamed Elhamdadi, Masahico Saito, Timothy Yeatman\\[2mm]
Department of Mathematics and Statistics\\
University of South Florida}

\date{\empty}

\maketitle

\begin{abstract}
We present a set of 26 finite quandles
that distinguish (up to reversal and mirror image) by number of colorings, all of the 2977 prime oriented knots with up to 12 crossings. We also show that 1058 of these knots can be distinguished from their mirror images by the number of colorings by quandles from a certain set of 23 finite quandles. We study the colorings of these 2977 knots by all of the 431 connected quandles of order at most 35 found by L. Vendramin. Among other things, we collect information about quandles that have the same number of colorings for all of the 2977 knots. For example, we prove that if $Q$ is a simple quandle of prime power order then $Q$ and the dual quandle $Q^*$ of $Q$ have the same number of colorings for all knots and conjecture that this holds for all Alexander quandles $Q$. We study a 
knot invariant based on a quandle homomorphism $f:Q_1\to Q_0$. We also apply the quandle colorings we have computed to obtain some new results for the bridge index, the Nakanishi index, the tunnel number, and the unknotting number.
In an appendix we discuss various properties of the quandles in Vendramin's list. Links to the data computed and various programs in C, GAP and Maple are provided.
\\

\noindent
{\it Key words}: Quandles, knot colorings, mirror images of knots, knot invariants.
\end{abstract}

\section{Introduction}\label{Introsec}

Sets with certain self-distributive operations called  {\it quandles}
have been studied since  the 1940s 
in various areas with different names. They have been studied,
for example, as  algebraic systems for symmetries \cite{Taka},
as quasi-groups~\cite{galkin}, and in relation to modules~\cite{XH,Sam}.
Typical examples of quandles arise from  conjugacy classes of groups
and from modules over the integral Laurent polynomial ring
$\Z [t, t^{-1}]$, called Alexander quandles. The {\it fundamental quandle} of a knot
was defined in a manner similar to the
fundamental group \cite{Joyce,Mat} of a knot, which made quandles an important
tool in knot theory.  The number of homomorphisms from the fundamental
quandle to a fixed finite quandle has an interpretation as colorings
of knot diagrams by quandle elements, and has been widely used as a
knot invariant. Algebraic homology theories for quandles 
were defined \cite{CJKLS,FRS1}, and
 investigated in \cite{LN,Moc,NP1,NP2,Nos}.
 Extensions of quandles by cocycles have been studied \cite{AG,CENS,Eis}, and
invariants derived thereof are applied to various properties of knots and knotted surfaces (see \cite{CKS}
and references therein).

Tables of small quandles have been made previously (e.g.,
\cite{CKS,EMR,Sam}).  Computations using GAP by L.  Vendramin
\cite{Vendramin} significantly expanded the list for connected quandles.  He found all
connected quandles of order up to $35$. There are $431$ of them. These
quandles may be found in the GAP package RIG \cite{rig}.  We refer to
these quandles as { \it RIG} quandles, and use the notation $C[n,i]$
for the $i$-th quandle of order $n$ in his list. As a matrix $C[n,i]$
is the transpose of the quandle matrix $SmallQuandle(n,i)$ in
\cite{rig}.

In this paper, we investigate to what extent the number of quandle
colorings of a knot by a finite quandle can distinguish the oriented
knots with at most 12 crossings from the knot table at KnotInfo
\cite{KI}.

In Section~\ref{prelimsec}, we provide definitions and conventions.
The coloring of knots by quandles and the relationship to symmetry of
knots are discussed in Section~\ref{distknots} and
Section~\ref{computation}. It turns out that {\it RIG} quandles do not
suffice to distinguish all of  the prime, oriented knots with at most 12
crossings.  To distinguish all of these knots we generated several
thousand conjugation quandles and generalized Alexander
quandles. Eventually we found a set of $26$ quandles that distinguish,
up to orientation and mirror image, all knots with up to $12$
crossings. These computations extend the results by Dionisio and Lopes
\cite{DL} for $10$ Alexander quandles and 249 prime knots with at most
10 crossings.

We write $m(K)$ for the mirror image of knot $K$ and $r(K)$ for $K$
with the orientation reversed. The 2977 knots given at KnotInfo
\cite{KI} are representatives up to mirrors and reverses
of the prime, oriented knots with at most 12 crossings (see \cite{symmetry_types}). 
It is known \cite{Man,Mat} that quandle
colorings do not distinguish $K$ from $rm(K)$ for any knot $K$.  Thus
it is of interest which finite quandles, if any, can distinguish $K$ from
$m(K)$.  This is only possible for a knot which is chiral or negative
amphicheiral. There are 1366 chiral or negative amphicheiral knots
among the 2977 knots with up to 12 crossings at KnotInfo \cite{KI}.
We have found a set of $23$ quandles that distinguish $K$ from $m(K)$
for 1058 knots out of these 1366 knots.
It remains an open question whether the remaining 308 of these knots can be so distinguished.

Many sets of RIG quandles share the same number of colorings for all knots in
the KnotInfo table.  We explore this phenomenon in
Section~\ref{similarities}. For example, we conjecture that a
connected Alexander quandle and its dual quandle always give the same
number of colorings for any knot, and we prove this for a special class of Alexander quandles.
 In Section~\ref{hom-sec}, we introduce an invariant of knots based on quandle homomorphisms and mention a relation to the quandle cocycle invariant.
 The computational results are applied to  other knot invariants in
Section~\ref{othersec}, such as the bridge index, the Nakanishi index, the tunnel number and the unknotting number.
Various properties of quandles are further discussed in
Appendix~\ref{propertysec}, and a summary of computational results is
given as to these properties for RIG quandles.  Descriptions  of quandles
used for distinguishing knots and their mirror images can be found in
Appendix \ref{tablesec}.
Files containing the quandle matrices, programs and the results of computations can be found at \cite{URL,BridgeIndex}.
Appendix~\ref{unk-sec} contains the list of 12 crossing knots discussed in Section~\ref{othersec}. 
 Appendix~\ref{computationsec} contains some information on  programs used for colorings.

\section{Preliminaries}\label{prelimsec}

We briefly review some definitions and examples of quandles. 
More details can be found, for example, in \cite{AG,CKS,FRS1}. 
 
A {\it quandle} $X$ is a 
set with a binary operation $(a, b) \mapsto a * b$
satisfying the following conditions.
\begin{eqnarray*}
& &\mbox{\rm (1) \ }   \mbox{\rm  For any $a \in X$,
$a* a =a$.} \label{axiom1} \\
& & \mbox{\rm (2) \ }\mbox{\rm For any $b,c \in X$, there is a unique $a \in X$ such that 
$ a*b=c$.} \label{axiom2} \\
& &\mbox{\rm (3) \ }  
\mbox{\rm For any $a,b,c \in X$, we have
$ (a*b)*c=(a*c)*(b*c). $} \label{axiom3} 
\end{eqnarray*}
 A {\it quandle homomorphism} between two quandles $X, Y$ is
 a map $f: X \rightarrow Y$ such that $f(a*_X b)=f(a) *_Y f(b) $, where
 $*_X$ and $*_Y$ 
 denote 
 the quandle operations of $X$ and $Y$, respectively.
 A {\it quandle isomorphism} is a bijective quandle homomorphism, and 
 two quandles are {\it isomorphic} if there is a quandle isomorphism 
 between them.

 Typical examples of quandles include the following. 
  \begin{itemize}
 \setlength{\itemsep}{-3pt}
\item
Any non-empty set $X$ with the operation $a*b=a$ for any $a,b \in X$ is
a quandle called a  {\it trivial} quandle.

\item
A conjugacy class $X$ of a group $G$ with the quandle operation $a*b=b^{-1} a b $.
We call this a {\it conjugation quandle}.

\item
A {\it generalized Alexander quandle} is defined  by  
a pair $(G, f)$ where 
$G$ is a  group and $f \in {\rm Aut}(G)$,
and the quandle operation is defined by 
$x*y=f(xy^{-1}) y $. 
If $G$ is abelian, this is called an {\it Alexander quandle}.

\item
A  {\it Galkin quandle}
is defined as follows. 
Let  $A$ be an abelian group, also regarded naturally as a $\Z$-module. 
Let $\mu: \Z_3  \rightarrow \Z$ ,  $\tau: \Z_3 \rightarrow A$ be functions.   
These functions $\mu$ and $\tau$ need not be homomorphisms.
Define a binary operation on $\Z_3 \times A$ by 
$$(x, a)*(y, b)=(2y-x, -a + \mu(x-y) b + \tau(x-y) ) \quad x, y \in \Z_3, \ a, b \in A. $$
Then for any abelian group $A$, 
the above operation $*$ defines a quandle structure on $\Z_3 \times A$ if 
$\mu(0)=2$, $\mu(1)=\mu(2)=-1$, and $\tau(0)=0$. 
Galkin gave this definition in \cite{galkin}, page $950$, for $A = \Z_p$,
and this definition was given in \cite{CEHSY}.

\item
A function $\phi: X \times X \rightarrow A$ for an abelian group $A$ is 
called a {\it quandle $2$-cocycle}  \cite{CJKLS} if it satisfies 
$$ \phi (x, y)- \phi(x,z)+ \phi(x*y, z) - \phi(x*z, y*z)=0$$
and $\phi(x,x)=0$ for any $x,y,z \in X$.
For a quandle $2$-cocycle $\phi$, 
$X \times A$ becomes a quandle 
by $(x, a) * (y, b)=(x*y, a+\phi(x,y))$ for $x, y \in X$, $a,b \in A$,
and it is called an {\it abelian extension} of $X$ by $A$, see  \cite{CENS}.

  \end{itemize}
  
  Let $X$ be a quandle.
  The {\it right translation}  ${\cal R}_a:X \rightarrow  X$, by $a \in X$, is defined
by ${\cal R}_a(x) = x*a$ for $x \in X$. Similarly the {\it left translation} ${\cal L}_a$
is defined by ${\cal L}_a(x) = a*x$. Then ${\cal R}_a$ is a permutation of $X$ by Axiom (2).
The subgroup of ${\rm Sym}(X)$ generated by the permutations ${\cal R}_a$, $a \in X$, is 
called the {\it inner automorphism group} of $X$,  and is 
denoted by ${\rm Inn}(X)$. 
A quandle is {\it connected} if ${\rm Inn}(X)$ acts transitively on $X$.
 The operation $\bar{*}$ on $X$ defined by $a\ \bar{*}\ b= {\cal R}_b^{-1} (a) $
is a quandle operation, and $(X,  \bar{*}) $ is called the {\it dual} quandle of $(X, *)$.
We also denote the dual of $X$ by $X^*$.
If  $X^*$ is isomorphic to  $X$, then $X$ is called {\it self-dual}.

\begin{figure}[ht]
\begin{minipage}[b]{0.5\linewidth}
\centering
\includegraphics[scale=.5]{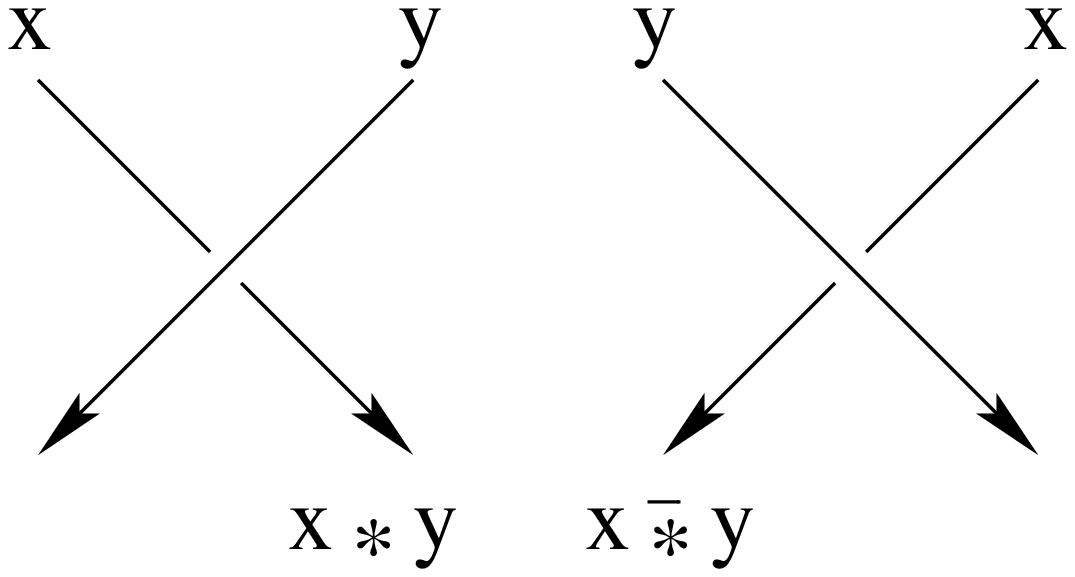}
\caption{The coloring rule at crossings}
\label{crossing}
\end{minipage}
\begin{minipage}[b]{0.5\linewidth}
\centering
\includegraphics[scale=.4]{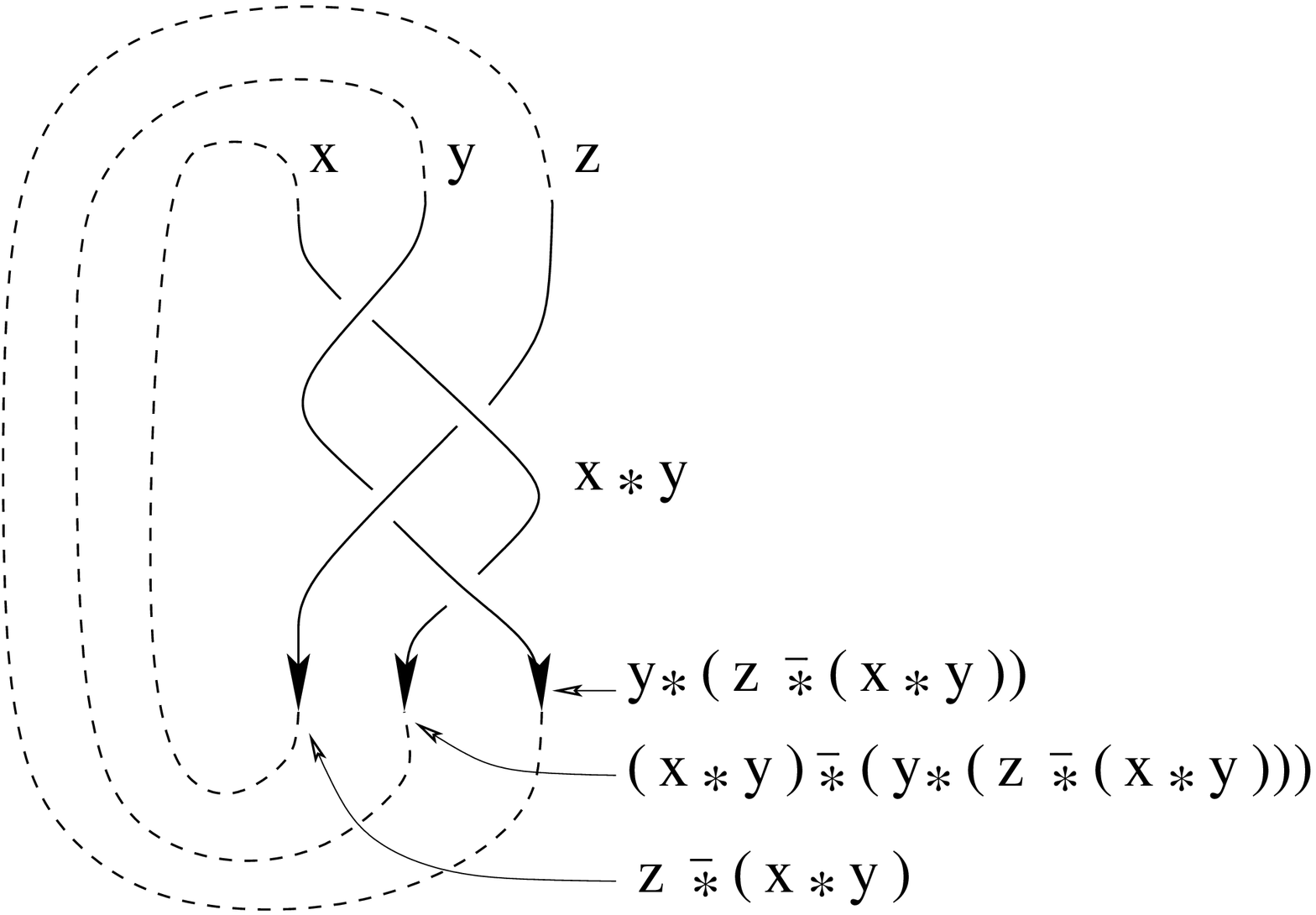}
\caption{Figure-eight   knot as the closure of the braid $\sigma_1 \sigma_2^{-1}  \sigma_1 \sigma_2^{-1}   $}
\label{figure8}
\end{minipage}
\end{figure}

 We write $K = K'$ to denote that there is an orientation preserving
 homeomorphism of $\Sp^3$ that takes $K$ to $K'$ preserving the
 orientations of $K$ and $K'$.  A {\it coloring} of an oriented knot
 diagram by a quandle $X$ is a map ${\cal C}: {\cal A} \rightarrow X$
 from the set of arcs ${\cal A}$ of the diagram to $X$ such that the
 image of the map satisfies the relation depicted in
 Figure~\ref{crossing} at each crossing.  More details can be found in
 \cite{CKS,Eis2}, for example.  A coloring that assigns the same
 element of $X$ to all the arcs is called trivial, otherwise
 non-trivial.  The number of colorings of a knot diagram by a finite
 quandle is known to be independent of the choice of a diagram, and
 hence is a knot invariant.
In Figure~\ref{figure8}, a coloring of a figure-eight knot is
indicated.  In the figure, the elements $x, y, z $ of a quandle are
assigned at the top arcs.  The coloring conditions are applied to
color other arcs as depicted.  Since the bottom arcs are connected to
the top arcs, we need the conditions
\begin{eqnarray*}
x &=& z \ \overline{*} \ (x * y) , \\
y &=& (x*y) \ \overline{*} \ ( y * ( z \ \overline{*} \ (x * y) ) ) ,\\
z &=& y * ( z \ \overline{*} \ (x * y) ) . 
\end{eqnarray*}
to obtain a coloring.
Each crossing corresponds to a standard generator $\sigma_i$ of the
$3$-strand braid group or its inverse, and the closure of $\sigma_1
\sigma_2^{-1} \sigma_1 \sigma_2^{-1} $ represents the figure-eight
knot.  The dotted lines indicate the closure, see \cite{Rolf} for more
details on braids.  We give the orientation downward for braids, see
the figure.  { \it For the rest of this paper all knots will be prime
  and oriented.}

The fundamental quandle is defined in a manner similar to the
fundamental group \cite{Joyce,Mat}.  A {\it presentation} of a quandle
is defined in a manner similar to groups as well, and a presentation
of the fundamental quandle is obtained from a knot diagram (see, for
example, \cite{FR}), by assigning generators to arcs of a knot
diagram, and relations corresponding to crossings.  The set of
colorings of a knot diagram $K$ by a quandle $X$, then, is in
one-to-one correspondence with the set of quandle homomorphisms from
the fundamental quandle of $K$ to $X$.

 If $Q$ is a quandle and $K$ is a
knot we denote by $\nkc{Q}{K}$ the number of non-trivial colorings of
$K$ by $Q$. The number of colorings including trivial colorings is 
$\kc{Q}{K} = \nkc{Q}{K} + |Q|$, where $|Q|$ is the order of $Q$. It
is well known that for each quandle Q, $\kc{Q}{K}$ is an invariant of
knots (\cite{Jozef}, for example). We say that the quandle $Q$ {\it
  distinguishes} knots $K$ and $K'$ if $\kc{Q}{K} \neq \kc{Q}{K'}$.

 Let $m: \Sp^3 \rightarrow \Sp^3$ be an orientation reversing
 homeomorphism.  For a knot $K$, $m(K)$ is the mirror image of $K$.
 Let $r(K)$ denote the same knot $K$ with its orientation reversed.
 We regard $m$ and $r$ as maps on equivalence classes of knots.  We
 consider the group $\mathcal{G} = \{ 1,r,m,rm \} $ acting on the set of all oriented knots.
  For each knot $K$ let $\mathcal{G}(K) =\{ K,r(K),m(K),rm(K)
 \}$ be the orbit of $K$ under the action of $\mathcal{G}$.

By a symmetry  we mean that a knot (type) $K$ remains unchanged under one of $r$, $m$, $rm$.  
As in the definition of {\it symmetry type} in \cite{KI} we say that  a knot $K$ is
\begin{itemize}
\setlength{\itemsep}{-3pt}
\item { {\it reversible} if the only symmetry it has is $K = r(K)$},
\item { {\it negative amphicheiral} if the only symmetry it has is $K = rm(K)$ },
\item {  {\it positive amphicheiral} if the only symmetry it has is $K = m(K)$},
\item {{\it chiral} if it has none of these symmetries},
\item {{\it fully amphicheiral} if it has all three symmetries, that is, $K =r(K)=m(K)= rm(K)$}.
\end{itemize}
 The symmetry type of each knot on at most 12 crossings is given at \cite{KI}. Thus each of the $2977$ knots $K$ given there represents as many as four knots $K, m(K), r(K)$ and $ rm(K)$.

It is known \cite{Man,Mat} that 
the fundamental quandles of $K$ and $K'$ are isomorphic if and only if 
$K=K'$ or $K=rm(K')$.

\section{Distinguishing the knots $K,m(K),r(K),rm(K)$}\label{distknots}

For the problem of distinguishing knots by quandle coloring  it suffices to consider
only connected quandles by the following lemma (see Ohtsuki
\cite{ohtsuki}, Section 5.2).

\begin{lemma}
The image of a quandle homomorphism of the fundamental quandle of any
knot in another quandle is connected.
In particular,  
given two knots $K_1$ and $K_2$, a quandle $Q$ of smallest order satisfying $\kc{Q}{K_1}  \neq \kc{Q}{K_2}$  is connected. 
\end{lemma}

\begin{lemma} \label{qk_equations} The following equations hold for all quandles $Q$ and all knots $K$:

\begin{list}{}{}
\setlength{\itemsep}{-3pt}

\item[$(1)$] { $\kc{Q}{K}=\kc{Q}{rm(K)}$},

\item[$(2)$] { $\kc{Q}{r(K))}= \kc{Q}{m(K))}$}, 

\item[$(3)$] {  $\kc{Q}{m(K)} = \kc{Q^*}{K}$}, 
\end{list}

where $Q^*$ in  $(3)$ denotes the dual quandle of $Q$. 

\end{lemma}

\noindent {\it Proof:} For (1) and (2) see  Kamada \cite{K}. Item (3) follows easily from the definitions. \qed

\medskip

\begin{corollary}\label{S_1}  For all  quandles $Q$ and for knots $K$ of any symmetry type except chiral or negative amphicheiral,
 $$\kc{Q}{rm(K)} = \kc{Q}{m(K)} = \kc{Q}{r(K)} = \kc{Q}{K}. $$ 
\end{corollary}

\begin{conjecture} \label{conj1}
  If $\mathcal{K}$ is any finite set of knots closed under the action %
  of $\mathcal{G}$, there exists a finite sequence of finite quandles
  $ \mathcal{S} = (Q_1, Q_2, \ldots, Q_{n}) $ such that the
  invariant $$C(K):=(\kc{Q_1}{K}, \kc{Q_2}{K}, \ldots,
  \kc{Q_{n}}{K})$$ satisfies for all $K,K' \in \mathcal{K}$:

  $$ C(K') = C(K) \text { if and only if } K' = K \text { or } K' = mr(K).$$
\end{conjecture}

Our computations verify the truth of the following weaker conjecture
where $\mathcal{K}$ is the set of all knots with at most 12
crossings. We also verify the truth of Conjecture \ref{conj1} for
smaller sets of knots (see Section \ref{computation}).

\begin{conjecture} \label{conj2}
  If $\mathcal{K}$ is any finite set of knots closed under the action %
  of $\mathcal{G}$, there exists a finite sequence of finite quandles
  $ \mathcal{S} = (Q_1, Q_2, \ldots, Q_{n}) $ such that the
  invariant $$C(K):=(\kc{Q_1}{K}, \kc{Q_2}{K}, \ldots,
  \kc{Q_{n}}{K})$$ satisfies for all $K,K' \in \mathcal{K}$:

  $$ C(K') = C(K) \text{ implies } \mathcal{G}(K') = \mathcal{G}(K)$$
\end{conjecture}
It is clear that if Conjecture \ref{conj1} is true then the following
conjecture is true.

\begin{conjecture} \label{conj3}

  There exists a sequence of finite quandles $ \mathcal{S} = (Q_1,
  Q_2, \ldots, Q_{n}, \ldots) $ such that the
  invariant $$C(K):=(\kc{Q_1}{K}, \kc{Q_2}{K}, \ldots,
  \kc{Q_{n}}{K},\ldots )$$ satisfies for all knots $K,K' $, 

  $$ C(K') = C(K) \text { if and only if } K' = K \text { or } K' = mr(K).$$
\end{conjecture}

This conjecture is similar to a conjecture made by Fenn and Rourke
\cite{FR}, page 385, for racks.

Let $\mathcal{K}$ be a finite set of knots closed under the action of the group $\mathcal{G}$.  
Let $\mathcal{R} = \{k_1,k_2, \ldots, k_n \}$ be a set of representatives of the orbits of $\mathcal{G}$ acting on $\mathcal{K}$. 
Let $ \mathcal{S} = (Q_1, Q_2, \ldots, Q_n)$  be a list of finite quandles and for any knot $K$ let
 $$C(K):=(\kc{Q_1}{K}, \kc{Q_2}{K}, \ldots, \kc{Q_{n}}{K}).$$
To simplify the following proof,  for knots $K$ and $K'$ we write $$ K \sim K' \text { if and only if } C(K) = C(K').$$
Thus $K \nsim K' \text { if and only if there exist } Q_i \in \mathcal{S} \text{ such that } \kc{Q_i}{K} \neq \kc{Q_i}{K'}$.
The  analysis of our computations is based on the following.

\begin{proposition} \label{th_conj1}

Staying with the notation above, suppose for all pairs of knots $k,k' \in \mathcal{R}$ the list $\mathcal{S}$ satisfies 
\begin{list} {}{}
\setlength{\itemsep}{-3pt}
\item[(A)]{$\text{if } k \neq k' \text{ then } k \nsim k'   $,}
\item[(B)]{$\text{if } k \neq k'   \text{ then } k \nsim m(k')$,  \text{ and } }
\item[(C)]{$\text{if } k \neq k'  \text{ then } m(k) \nsim m(k') $.}
\end{list}
Then if $K,K' \in \mathcal{K}$ and $K \sim K'$ we have $ \mathcal{G}(K) = \mathcal{G}(K')$.  If in addition 
 $\mathcal{S}$ satisfies
\begin{list}{}{}
\item[(D)] {$m(k) \nsim k \text{ when $k$ is chiral or negative amphicheiral}$,}
\end{list}
then if  $K,K' \in \mathcal{K}$ and $K \sim K'$ we have  $K' = K$ or $K' = rm(K)$.

\end{proposition}
\noindent {\it Proof:} Assume first that (A), (B), (C) hold and that $K,K' \in \mathcal{K}$  and $K \sim K'$.  We must prove that $\mathcal{G}(K) = \mathcal{G}(K')$. Let $k,k'$ be elements of $\mathcal{R}$ such that $$K \in \mathcal{G}(k) \text{ and  } K' \in \mathcal{G}(k').$$ 
Note that 
$\mathcal{G}(k) = \{ k,rm(k) \} \cup \{ r(k),m(k) \} $ and 
$ \mathcal{G}(k') = \{ k',rm(k') \} \cup \{ r(k'),m(k') \}. $ 
By Lemma \ref{qk_equations} (1) and  (2)  we have $k \sim rm(k)$ and $r(k) \sim m(k)$,
and similarly for $k$ replaced by $k'$.
We  consider four cases and show that each case leads to $\mathcal{G}(k)=\mathcal{G}(k')$. 

 \begin{list}{}{}
\setlength{\itemsep}{-3pt}
\item[Case(i)] { $K \in \{k,rm(k) \} \text{ and } K' \in \{k',rm(k' ) \}$. 
Then we have $k \sim K \sim K' \sim k'$. This implies by (A) that $k =  k'$, and hence 
$\mathcal{G}(k)=\mathcal{G}(k')$.}
\item[Case(ii)] { $K \in \{k,rm(k) \}  \text{ and } K' \in \{m(k'),r(k' ) \}$.  Then $ k \sim K \sim K' \sim m(k')$ which implies by (B) that $k = k'$, hence $\mathcal{G}(k)=\mathcal{G}(k')$.}
\item[Case(iii)] { $K \in \{m(k),r(k) \}  \text{ and } K' \in \{k',rm(k' ) \}$. This is similar to  Case(ii).}
\item[Case(iv)] { $K \in \{m(k),r(k) \}  \text{ and } K' \in \{m(k'),r(k' ) \}$.  Then $m(k) \sim K \sim K' \sim m(k')$,  and  by (C) we have
$k =  k'$, hence $\mathcal{G}(k)=\mathcal{G}(k')$. }
\end{list}

This proves the first part of the proposition. 
Now assume that (A), (B), (C), (D) hold and $K \sim K'$. Let $k$ and $k'$ be as above.  By
what we just proved we have $\mathcal{G}(k)=\mathcal{G}(k')$ and so we have $K,K' \in \mathcal{G}(k)$.
If $K = K'$ there is nothing to prove so we can assume that one of the following 4 cases holds and that $K \neq K'$.
 \begin{list}{}{}
\setlength{\itemsep}{-3pt}
\item[Case(1)] {$k$ is reversible or positive amphicheiral.  Hence $\mathcal{G}(k)= \{ k, rm(k) \}$. In this case we  may take $K = k$ and $K' = rm(k)$ so $K' = rm(K)$.}
     
\item[Case(2)] {$k$ is negative amphicheiral. Hence $\mathcal{G}(k)= \{ k, m(k) \}$.  In  this case by (D) we have that $k \nsim m(k)$ so $K \neq K'$,  $K \sim K'$ and
$K,K' \in \mathcal{G}(k)$ is impossible.}
\item[Case(3)] {$k$ is chiral. Hence $\mathcal{G}(k)= \{ k,r(k), m(k),rm(k) \}$.  Of the four elements of $\mathcal{G}(k)$ there are 6 possible pairs that might be $\{K,K'\}$.
Note that the pair $\{k,m(k)\}$ is ruled out by (D). Since we always have $k \sim rm(k)$ and $r(k) \sim m(k)$,
the pairs $\{k,r(k) \}$, $\{ rm(k), m(k) \}$, $\{ rm(k),r(k) \}$, are also  ruled out. This leaves the only possibilities to be $\{ k, rm(k) \}$ and $\{r(k),rm(k)\}$. If $r(k) \sim rm(k)$ then $m(k) \sim r(k) \sim rm(k) \sim k$ which is a contradiction to (D), so we are left with $\{k,rm(k) \}$ as the only possibility for $\{K,K' \}$, as desired.}
\item[Case(4)] {$k$ is fully amphicheiral. Hence $\mathcal{G}(k)= \{ k \}$. This is impossible since there is only a single knot in $\mathcal{G}(k)$.}
\end{list}

This completes the proof.  \qed

\section{Computational support for Conjectures \ref{conj1}  and \ref{conj2}} \label{computation}

Let $\mathcal{Q}_{26} = ( Q_1, \ldots , Q_{26} )$ be the list of 26 quandles in Table~\ref{table-Q} in Appendix B. 
Matrices for these 26 quandles may be  found  at  \cite{URL}, 
where they are denoted by $Q[i], i = 1, \ldots, 26$.  The list
$\mathcal{Q}_{26}$ is closed under taking dual quandles. This is to
simplify computing the number of colorings of the mirror image of a
knot. Included in the file is the list {\it dual} satisfying
$dual[i]=j$ if and only if $Q_i^* $ is isomorphic to $ Q_j.$

If $i \neq 15$ each quandle $Q_i$ is faithful and hence isomorphic to
a conjugation quandle on some conjugacy class in its inner
automorphism group Inn($Q_i$). The quandle $Q_{15}$ is not faithful
and is not a conjugation quandle [Vendramin, personal communication].
We note that some
non-faithful quandles are conjugation quandles, so this is not
trivial. 
However, $Q_{15}$ is a generalized Alexander quandle on GAP's
SmallGroup$(32, 50) = (C2 \times Q8 ) :C2$. 
 The quandles in ${\mathcal Q}_{26}$ of order less than 36 are RIG quandles and the
remaining 11 were found by searching for conjugation quandles on
conjugacy classes of finite groups using GAP.  In retrospect all of
these 26 quandles could have been found by searching for conjugation
quandles and generalized Alexander quandles.

We denote the knots in KnotInfo up to 12 crossings
by $$\mathcal{R}_{12}=\{ K[i] : i = 1, \ldots, 2977 \}.$$ In our
notation, the knot $K[i]$ represents the knot with {\it name rank} $=
i+1$ in \cite{KI}.  In particular, $K[1]$ is the trefoil, which has
name rank $2$ in \cite{KI}.

The $26 \times 2977$ matrix $M$ found at
\cite{URL} has the property that $$M_{i,j} = \nkc{Q_i}{K[j]}.$$ See
Appendix~\ref{computationsec} for comments on the computation of
$\nkc{Q}{K}$.

Since $\mathcal{R}_{12}$ is a set of representatives of the set
$\mathcal{K}_{12}$ of (isotopy classes of prime, oriented) knots with
at most $12$ crossings, including all symmetry types, to verify that
Conjecture \ref{conj1} holds when $\mathcal{K} = \mathcal{K}_{12}$ and
the list of quandles $\mathcal{S} = \mathcal{Q}_{26}$, it suffices, by
Proposition \ref{th_conj1}, to verify that the matrix $M$ and the list
$dual$ has the following properties for $i,j \in \{1,2, \ldots,
2977\}$ and $k \in \{1,2,\ldots,26 \}$.

\begin{list}{}{}
\setlength{\itemsep}{-3pt}
\item[$(1)$] {For all $i \neq j$ there exists  $k$ such that $M_{k,i} \neq M_{k,j }$.}
\item[$(2)$]  {For all $i \neq j$ there exists  $k$ such that $M_{k,i} \neq M_{s,j }$, where $s = dual[k]$.}
\item[$(3)$]  {For all $i \neq j$ there exists  $k$ such that $M_{s,i} \neq M_{s,j }$, where $s = dual[k]$.}
\end{list}

These conditions are easy to verify directly from the matrix $M$ and
the list $dual$. This shows that the set $\mathcal{Q}_{26}$
distinguishes the knots in $\mathcal{K}_{12}$ up to mirror images and
orientation.

\begin{remark}
{\rm 
Note that the unknot is not included among the knots K[i], $i = 1,
\ldots , 2977$. Since the unknot has only trivial colorings, to show
that we can distinguish it from all of the 2977 knots it suffices to
note that for every knot $K[j]$ there is a quandle $Q_i$ such that
$\nkc{Q_i}{K[j]} \neq 0$, that is, no column of the matrix M contains
only $0$s.  This is confirmed from the matrix $M$. 
}
\end{remark}

As pointed out above the only knots we can hope to distinguish from
their mirror images (equivalently from their reversals) are chiral
knots and negative amphicheiral knots. Among 2977 knots in KnotInfo \cite{KI}
up to 12 crossings, 1319 knots are chiral and 47 are negative
amphicheiral. Denote this set of 1366 knots by $\mathcal{K}^c$.  We
were not able to find quandles to distinguish all of these 1366 knots
from their mirror images.  Thus we cannot verify
Conjecture~\ref{conj2} when $\mathcal{K} = \mathcal{K}_{12}$.
However, were able to find quandles that distinguish 1058 of the knots
in $\mathcal{K}^c$ from their mirror images.  The quandles that
distinguish these 1058 knots are the quandles $T_1, T_2,
\ldots,T_{23}$ in Table~\ref{table-T}.
Matrices for the quandles can be found in \cite{URL}.  
A list of 1058 pairs $(i,n)$ where $K[n]$  is one of 1058
 knots in $\mathcal{K}^c$ and $T_i$ is a quandle that distinguished
 knot $K[n]$ from $m(K[n])$ can be found in \cite{URL}. 

With the exception of quandles $T_2$ and $T_3$ the quandles $T_i$ are
conjugation quandles on conjugation classes in the group
$Inn(T_i)$. $T_2$ and $T_3$ are generalized Alexander quandles, both
on the group $PSL(2,8)$. We note that $Q_{17} = T_{22}$ and $Q_{24} =
T_{23}$. Otherwise the quandles in Table~\ref{table-T} are different
from those in Table~\ref{table-Q}

In Table~\ref{table-Q} and Table~\ref{table-T} the structure
descriptions are given by GAP's {\it StructureDescription} function. In particular, $A\times B$ is the
direct product, $N:H$ is a semidirect product and $A.B$ is a non-split
extension. The group ID, $SmallGroup(n,i)$ is the $i$th group of order
$n$ in the GAP Small Groups library for $n$ at most 2000. If a group
has order $n > 2000$ we write $SmallGroup(n,0)$. Nevertheless for all
of the  inner automorphism groups we have computed, GAP does give a Structure Description.

\section{Similarity of quandles over knot colorings} \label{similarities}

 For two quandles $Q_1$ and $Q_2$, and a set
${\mathcal K}$ of knots, we write 
$Q_1 \approx_{\mathcal K} Q_2$ if $\kc{Q_1}{K} = \kc{Q_2}{K} $ for all $K \in {\mathcal K} $.  

This equivalence relation was considered in \cite{CEHSY}.  We omit
the subscript ${\mathcal K}$ when it is the set of all knots.  The
following is immediate from the definitions and a property of the
number of colorings.

\begin{lemma}\label{prod-lem}
If $Q_1 \approx Q_1'$ and $Q_2 \approx Q_2'$ 
 then 
$( Q_1 \times Q_2)  \approx (Q_1' \times  Q_2')$. 

\end{lemma}

Let ${\mathcal K}$ be the set of all  $2977$ knots in the table in
KnotInfo \cite{KI} up to $12$ crossings.
We observe the following
for $\approx_{\mathcal K}$ and for the 431 RIG quandles.

\begin{itemize}
\setlength{\itemsep}{-3pt}
\item There are $151$ classes of $\approx_{\mathcal K}$ consisting of more than one quandle.

\item
Among these classes, $145$ classes consist of two  quandles. 

\item
Among $145$ classes of pairs, all but $35$ classes consist of a pair of dual Alexander quandles.

\item
Among $35$ classes of pairs,  one class is  a pair of  self-dual Alexander quandles, 
some are non-Alexander dual quandles, 
and some are non-Alexander self-dual quandles.

\end{itemize}

\begin{conjecture}{\rm 
For every connected Alexander quandle $Q$, $Q \approx Q^*$. 
}
\end{conjecture}

 We prove the above conjecture  for simple Alexander quandles  of prime power order.  (Following \cite{AG} we say that a quandle $Q$ is {\it simple} if it is not trivial and whenever $f:Q \to Q'$ is an epimorphism then $f$ is an isomorphism or  $|Q'| = 1$. This implies that $|Q| > 2$ since the quandles of orders 1 and 2 are trivial.) The following proposition 
 gives  alternative characterizations of these quandles.

\begin{proposition}
\label{ItalianQuandles}
The following statements are equivalent for a quandle $Q$ and a prime $p$.
\begin{list}{}{}
\setlength{\itemsep}{-3pt}

\item[(1)]{$Q$ is a simple quandle with $p^k$ elements.}
\item[(2)]{$Q$ is isomorphic to an Alexander quandle on $\Z_p[t, t^{-1}]/(h(t))$ where $h(t) \in \Z_p[t] $ is irreducible, of degree $k$ and different from  $t$ and $t-1$  and the quandle operation is defined by 
$x*y = \overline{t}x + (1-\overline{t})y $ where $\overline{t}$ is the coset $t+(h(t))$.}
\item[(3)]{$Q$ is isomorphic to an Alexander quandle on the vector space ${(\F_p)}^k$ where $T$ is the companion matrix of a monic irreducible polynomial $h(t) \in \Z_p[t] $ of degree $k$, different from $t$ and $t-1$ and the quandle operation is given by $x*y = Tx + (1-T)y$.}
\item[(4)]{$Q$ is isomorphic to a quandle defined on the finite field $\F_{q}$ where $q = p^k$,  $a \in \F_q$, $a \neq 1$, $a$ generates $\F_q$, and the quandle operation is given by $x*y=ax+(1-a)y$ for $x,y \in \F_q$.}
\end{list}

\end{proposition}
{\it Proof.\/} The equivalence of (1), (3) and (4) is proved in  \cite{AG} Section 3.
The equivalence of (2) and (4) is proved in \cite{BF}. 
$\Box$

\begin{proposition}\label{dualAlex}
If $Q$ is a simple  quandle of  prime power order,  
then $Q \approx Q^*$.
\end{proposition}
{\it Proof.\/}
By Proposition~\ref{ItalianQuandles}, 
the quandle $Q$ is equal to $(\F_q, *)$ where 
$\F_q =\Z_p[t, t^{-1}]/(h(t))$, $q=p^k$, 
 $h(t) \in \Z_p[t] $ is irreducible of degree $k$. 
 The quandle operation is defined by 
$x*y = \overline{t}x + (1-\overline{t})y $ where $\overline{t}$ is the coset $t+(h(t))$.
{}Let  $S \in M_n(\F_q )$ denote the reduction modulo $(p, h(t))$ of a Seifert matrix  of $K$, and 
let $N = \overline{t} S - S^T $, which is denoted by $N(1,p, h(t) )$ in \cite{BF}.
Although $p$  is restricted to odd primes  in \cite{BF}, we claim that their results hold for $p=2$
also.  The authors of \cite{BF} agree (personal communication).

Let 
$\varphi  : \F_q^ n \rightarrow \F_q^n  $ be the linear map that $N$ represents. 
Then in \cite{BF} it was shown that the set of colorings $\skc{Q}{K}$ is in one-to-one correspondence 
with the module $ \F_q \oplus {\rm ker}(\varphi) $. 

The dual quandle $Q^*$ is equal to $(\F_q, \overline{*})$ where 
$a \ \overline{*}\ b = {\overline{t}}^{-1} a + (1 - {\overline{t}}^{-1})b$.
Let $N^*=\overline{t}^{-1} S - S^T $, and $\varphi^*: \F_q^n \rightarrow \F_q^n$ be the corresponding linear map.
Then the number of colorings by $Q^*$ is the cardinality of $ \F_q \oplus {\rm ker}(\varphi^*) $.
It remains to show that the dimensions of ${\rm ker}(\varphi) $ and  ${\rm ker}(\varphi^*) $ are the same.
One computes  
$$(-\overline{t}^{-1} ) N^T = - \overline{t}^{-1} ( \overline{t}S^T -S)=-S^T + \overline{t}^{-1} S = N^*. $$
Since ${\rm rank}(N^T)={\rm rank}(N)$, we have  ${\rm rank}(N)={\rm rank}(N^*)$ and
hence the dimension of the kernel of $\varphi^*$ is equal to the dimension of the kernel of $\varphi$, so 
the number of  colorings is the same for $Q$ and $Q^*$.
$\Box$

This proof is similar to the proof of the fact that the Alexander polynomial is symmetric, $\Delta(t^{-1}) = \Delta(t)$ modulo a power of $t$. 

From Proposition \ref{dualAlex} and Lemma \ref{prod-lem} we obtain the following corollary.

\begin{corollary}
If a quandle $Q$ is isomorphic to a product of simple quandles of prime power order then
 $Q \approx Q^*$.
\end{corollary}

There are 273 connected Alexander quandles of order at most 35. Only 45 of these quandles are not  a product of simple quandles of prime power order.

\begin{remark}
{\rm All connected quandles of prime order clearly satisfy the conditions of Proposition \ref{ItalianQuandles}.
 Among RIG quandles, other than quandles of prime orders, the following quandles 
fall into this family.
\begin{eqnarray*}
 & & 
 C[ 4, 1 ],
C[ 8, 2 ],
C[ 8, 3 ],
C[ 9, 3 ],
C[ 9, 7 ],
C[ 9, 8 ] ,
C[ 16, 3 ],
C[ 16, 8 ],
C[ 16, 9 ],
C[ 25, 3 ],
\\
& & 
C[ 25, 11 ],
C[ 25, 12 ],
C[ 25, 13 ],
C[ 25, 19 ],
C[ 25, 20 ],
C[ 25, 31 ],
C[ 25, 32 ],
C[ 25, 33 ],
C[ 25, 34 ],
\\
& & 
C[ 27, 31 ],
C[ 27, 32 ],
C[ 27, 33 ],
C[ 27, 34 ],
C[ 27, 62 ],
C[ 27, 63 ],
C[ 27, 64 ],
C[ 27, 65 ],
C[ 32, 10 ],
\\
& & 
C[ 32, 11 ],
C[ 32, 12 ],
C[ 32, 13 ],
C[ 32, 14 ],
C[ 32, 15 ].
 \end{eqnarray*}

}
\end{remark}

We now  explain to the extent possible the 6 non-trivial classes of $\approx_{\mathcal K}$ that 
have more than two elements. 
Among these 6 classes, the first three classes listed below form equivalence classes under $\approx$ by Lemma~\ref{prod-lem} and Proposition~\ref{dualAlex}.
In all cases unless otherwise stated we have no explanation for an equivalence class.

\medskip

$\bullet$ \ $\{ C[25,6],C[25,7],C[25,8] \} $:  This situation is explained as follows: $C[25, 6]=C[5, 2] \times C[5, 2]$,
$C[25, 7]=C[5, 3] \times C[5, 3]$,
$C[25, 6]=C[5, 2] \times C[5, 3]$,
where 
$C[5, 2] =\Z_5 [t]/ (t+2)$ is dual to 
$C[5, 3] =\Z_5 [t]/ (t+3)$.

\medskip

$\bullet$ \ $\{ C[35,8],C[35,9],C[35,10],C[35,11] \} $:    This situation is explained as follows: 
$C[35, 8]=C[5, 2] \times C[7, 4]$,
$C[35, 9]=C[5, 2] \times C[7,5]$,
$C[35, 10]=C[5, 3] \times C[7, 4]$,
$C[35, 11]=C[5, 3] \times C[7,5]$,
where 
$C[7, 4] =\Z_7 [t]/ (t+2)$ is dual to 
$C[7, 5] =\Z_7 [t]/ (t+4)$.

\medskip

$\bullet$ \  $\{ C[35,12],C[35,13],C[35,14],C[35,15] \}$:  
This situation is explained as follows: 
$C[35, 12]=C[5, 3] \times C[7, 3]$,
$C[35, 13]=C[5, 3] \times C[7,2]$,
$C[35, 14]=C[5, 2] \times C[7, 3]$,
$C[35, 15]=C[5, 2] \times C[7,2]$,
where 
$C[7, 2] =\Z_7 [t]/ (t+3)$ is dual to 
$C[7, 3] =\Z_7 [t]/ (t+5)$.

\medskip

$\bullet$ \  $ \{ C[30,7],C[30,8],C[30,9],C[30,10] \} $:  This situations is explained as follows:
 $C[30,7]=C[5, 2] \times C[6,1]$,
$C[30,8]=C[5, 3] \times C[6,1]$,
$C[30,9]=C[5, 2] \times C[6,2]$,
$C[30,10]=C[5, 3] \times C[6,2]$,
where 
$C[5, 2] =\Z_5 [t]/ (t+2)$ is dual to 
$C[5, 3] =\Z_5 [t]/ (t+3)$. As mentioned before it is conjectured \cite{CEHSY} that $C[6,1] \approx C[6,2]$.
So we conjecture that this forms a class for $\approx$.

\medskip

$\bullet$ \  $\{ C[24,5],C[24,6],C[24,16],C[24,17] \}$:
These quandles are abelian extensions of $C[12,8]$, $C[12,9]$, $C[12,8]$,
$C[12,8]$, respectively, and none is Alexander.  The quandles
$C[12,8]$, $C[12,9]$ have $C[6,1]$, $C[6,2]$ as subquandles,
respectively, and have $C[6,1]$ as an epimorphic image. These facts alone, however, do not seem to 
imply that they constitute a class for $\approx$.
We conjecture
that this forms a class for $\approx$

\medskip

$\bullet$ \ $\{ C[24,29],C[24,30],C[24,31] \}$: These three quandles are self-dual Galkin quandles. This class was reported in \cite{CEHSY}.

\medskip

The remaining classes consist of pairs, that are not dual Alexander quandles, and we categorize them as follows. We conjecture that  they 
form equivalence classes under $\approx$.

There is one class with a pair of Alexander self-dual quandles:
$C[27,17]$  and $C[27,22]$, 
both of which  are  Alexander quandles on  the abelian group $\Z_3 \times \Z_9$.

Each of the following classes is a pair of  non-Alexander quandles that are duals of each other:
\begin{eqnarray*}
& & \{C[27,35],C[27,36] \},  \{C[27,37],C[27,38] \},  \{C[27,41],C[27,42] \} ,\\ 
& & \{C[27,43],C[27,46] \},  \{ C[27,44],C[27,45] \} ,  \{C[27,56],C[27,57] \}, \\
& & \{ C[27,58],C[27,59] \}, \{ C[28,11],C[28,12]\},  \{C[30,17],C[30,18] \},  \\
& & \{C[32,7],C[32,8]\}.
\end{eqnarray*}

The remaining classes are  pairs of self-dual  non-Alexander quandles. 
In \cite{CEHSY}, the relation $\approx$ was studied for a family called Galkin quandles.
Galkin quandles are known to be self-dual. 

The following are two element $\approx_{\mathcal K}$-equivalence classes of (non-Alexander) Galkin quandles that were reported in \cite{CEHSY}:
\begin{eqnarray*}
& &  \{ C[6,1],  C[6,2] \},
\{ C[12,5] , C[12,6] \}, 
\{ C[12,8],  C[12,9]  \}, \\
& & \{C[18,1] ,C[18,4] \}, 
 \{C[18,5] ,C[18,8] \}, 
\{C[24,27] ,C[24,28]\}, \\
 & &\{ C[24,38] ,C[24,39] \}, 
\{ C[30,12] ,C[30,14] \},
\{C[30,13] , C[30,15] \}. 
\end{eqnarray*}

The following classes are pairs of  self-dual  non-Alexander, non-Galkin quandles: 
\begin{eqnarray*}
& & \{C[18,2],C[18,3] \}, \{C[18,6],C[18,10] \},  \{ C[18,7],C[18,9] \}, \{C[20,5],C[20,6] \},\\
& & \{C[20,9],C[20,10] \} ,  \{ C[24,3],C[24,4] \},  \{ C[24,15],C[24,18] \},  \{ C[24,22],C[24,23] \}, \\
& & \{ C[24,34],C[24,36] \}, \{ C[24,35],C[24,37] \}, \{C[24,41],C[24,42] \}, \\ 
& & \{C[30,2],C[30,6] \}, \{C[30,23],C[30,24] \}, \{C[32,5],C[32,6] \}.
\end{eqnarray*}

\section{Quandle homomorphisms and knot colorings}\label{hom-sec}

We write $\skc{Q}{K}$ to denote the set of all colorings of the knot
$K$ by the quandle $Q$.

Let $f: Q_1 \rightarrow Q_0$ be a quandle homomorphism, and $K$ be a
knot.  Note that for any coloring $y \in \skc {Q_1}{K}$ regarded as a
quandle homomorphism $Q(K) \rightarrow Q_1$, the composition $x=f
\circ y :Q(K) \rightarrow Q_0$ is a coloring in $\skc{Q_0}{K}$.  In
this case we say that $x$ lifts to $y$.  This correspondence defines a
map $f_{\sharp}: \skc{Q_1}{K} \rightarrow \skc{Q_0}{K}$.

For $x \in \skc{Q_0}{K}$, denote by $f_{\sharp}^{-1}(x) \subset
\skc{Q_1}{K}$ the set of colorings $y \in \skc{Q_1}{K}$ such that
$f_{\sharp}(y)=x$. (This set may be empty.)  A pair of colorings $x
\in \skc{Q_0}{K}$ and $y \in \skc{Q_1}{K}$ such that $x=f \circ y $
for a knot diagram determines another such a pair after each
Reidemeister move, and hence for each $x \in \skc{Q_0}{K}$, the
cardinality $| f_{\sharp}^{-1}(x) |$ does not depend on the choice of
diagrams.  Thus we obtain the following knot invariant.

\begin{definition}
{\rm
For a given quandle homomorphism $f: Q_1 \rightarrow Q_0$, 
  we denote by ${\rm Col}_{f} (K) $ the multiset $ \{
  |f_{\sharp}^{-1}(x) | \ | \ x \in \skc{Q_0}{K} \}$, that does not
  depend on the choice of a diagram of a knot.
}
\end{definition}

We use the notation $[a,k]$ for $k$ copies of the element $a$. 
For example, $\{ [0, 3], [2,4] \}$ represents
$\{ 0,0,0,2,2,2,2\}$.
We observe the following.

\begin{itemize}
\item
If ${\rm Col}_{f} (K) =\{ [h_1, k_1], \ldots, [h_n, k_n] \}$,
then $\kc{Q_0}{K}= k_1 + \cdots + k_n$ and
 $\kc{Q_1}{K}= h_1 k_1 + \cdots + h_n k_n$. 

\item
Suppose $Q_1$ is connected. 
If $\nkc{Q_0}{K}=0$ (i.e.\ $K$ is only trivially colored by $Q_0$)
and $f$ is an epimorphism,
then ${\rm Col}_{f} (K) =\{ [h,k] \}$,
where $k=|Q_0|$ and $h=\kc{Q_1}{K}/|Q_0|$. 

%

\item
Suppose $Q_1$ is connected. 
If $\nkc{Q_1}{K}=0$, 
then ${\rm Col}_{f} (K) =\{ [0,\ell] , [h, k] \}$,
where $\ell=\nkc{Q_0}{K}$, $k=|Q_0|$, and $h=|Q_1|/|Q_0|$.


\item
If $f: Q_1 \rightarrow Q_0 $ has image  $Q_0' $, 
let $f': Q_1 \rightarrow Q_0' $ be the same map with the codomain $Q_0'$, 
then ${\rm Col}_{f} (K) ={\rm Col}_{f' } (K) \cup \{ [0, h] \}$ where 
$h$ is the number of colorings by $Q_0$ such that some colors are not from $Q_0'$.

\item 
If $f: Q_1=Q_0 \times Q_0' \rightarrow Q_0$ is the projection from a product quandle, 
then ${\rm Col}_{f} (K) =\{ [h, k] \}$, where $k=\kc{Q_0}{K}$ and $h=\kc{Q_0'}{K}$.

\end{itemize}

\begin{example}
{\rm 
We computed the invariant 
${\rm Col}_{f}  (K) $ 
for $i=1,2$ and for $2977$ knots up to $12$ crossings,
where $f$ is the unique  (up to isomorphisms of $C[6,i]$ and $C[3,1]$)
epimorphism from $C[6,i]$ to $C[3,1]$.
The values of the invariants for $i=1,2$ are the same for each $K$ up to $12$ crossings. 

%
%

For 2977 knots, the types of multisets are
\begin{eqnarray*}
& & \{[2, 3]\}, \ 
 \{[2, 3], [4, 6]\}, \ 
 \{[2, 3], [4, 24]\}, \
 \{[2, 3], [8, 6]\},  \
 \{[2, 3], [8, 24]\}, \\
& &  \{[2, 3], [4, 6], [8, 18]\}, \
 \{[2, 3], [4, 12], [8, 12]\} ,\
 \{[2, 3], [4, 18], [8, 6]\}, \\
& & \{[2, 3], [8, 18], [16, 6]\} , \
 \{[2, 3], [4, 36], [8, 24], [16, 18]\} .
\end{eqnarray*}

}
\end{example}

\begin{example}
{\rm 
There is an epimorphism $f: C[8,1] \rightarrow C[4,1]$ 
such that there are $12$ non-trivial colorings of 
the trefoil $K=3_1$  by $C[4,1]$, none of which lifts to a coloring 
by $C[8,1]$.
Hence for this $f:C[8,1] \rightarrow C[4,1]$ and for $K=3_1$,
the invariant  ${\rm Col}_{f: Q_1\rightarrow Q_0} (K)$ contains $12$ copies of $0$s.
There are $4$ and $8$  trivial colorings of $K$ by $C[4,1]$ and $C[8,1]$, respectively, 
and every trivial coloring by $C[4,1]$ has two lifts. 
Therefore we have 
 ${\rm Col}_{f} (K)= \{ [0, 12], [2,4] \}$.

We note that $C[8,1]$ is not Alexander but is an abelian extension of $C[4,1]=\Z_2[t]/(t^2+t+1)$.
Hence this computation is equivalent to quandle cocycle invariant.
All possible patterns over all of  the first $1000$ knots are
$$ \{[2, 4]\} , \ 
\{[2, 16]\} , \ 
 \{[2, 64]\} , \ 
 \{[0, 12], [2, 4]\} , \
 \{[0, 24], [2, 40]\} , \ 
\{[0, 48], [2, 16]\} .
$$

}
\end{example}

\begin{example}
{\rm 
For $f: C[9,1]=R_9 \rightarrow C[3,1]=R_3$, the patterns over all of  the first $1000$ knots  are
$$\{[3, 3]\} ,\ 
\{[9, 9]\} ,\
\{[27, 27]\} ,\
\{[0, 6], [9, 3]\} ,\
\{[0, 18], [27, 9]\} ,\
\{[0, 24], [27, 3]\} .$$

}
\end{example}

\begin{example}
{\rm 

There is an epimorphism $f: C[18,5] \rightarrow C[6,2] $.
The quandle $C[18,5]$ is  not Alexander,  not a Kei, not Latin, self-dual, and faithful.
It is not an abelian extension of $C[6,2]$. It has $3$ subquandles, $C[3,1]$, a
trivial quandle of $2$ elements, a $6$-element non-connected quandle,
and does not have $C[6,1]$ nor $C[6,2]$ as subquandles.
The RIG quandles onto which it has epimorphisms are $C[3, 1]$, $C[6, 2]$, and $C[9, 6]$. 

For this $f: C[18,5] \rightarrow C[6,2] $ based on  the first 100 knots, the patterns 
of the invariant are 
$$
\{[3, 6]\} ,\
\{[9, 30]\}  ,\
\{[9, 54]\}  ,\
\{[27, 102]\}  ,\
\{[0, 96], [27, 6]\}  ,\
\{[0, 120], [27, 6]\}.
$$

}
\end{example}

\begin{remark}
{\rm

For all examples of epimorphisms  $f: Q_1 \rightarrow Q_0$ we computed,
the invariant 
 $\kc{f}{K}$ is determined by $(\kc{Q_1}{K}, \kc{Q_0}{K} )$. 
Specifically, for these epimorphisms and all pairs of knots $K$ and $K'$ 
with at most $12$ crossings, if $(\kc{Q_1}{K}, \kc{Q_0}{K} ) = (\kc{Q_1}{K'}, \kc{Q_0}{K'} ) $
then  $\kc{f}{K} = \kc{f}{K'} $. 
This holds for epimorphisms 
$C[6,1] \rightarrow C[3,1]$, 
$C[8,1] \rightarrow C[4,1]$, and
$C[9,1] \rightarrow C[3,1]$.
It also holds for 
$C[18,5] \rightarrow C[6,2]$, 
for many  knots (not necessarily all $2977$ knots) we computed.

Furthermore, for $C[6,1] \rightarrow C[3,1]$, the invariant is determined by $\kc{C[6,1]}{K}$ alone for many knots, 
while for $C[8,1] \rightarrow C[4,1]$ and
$C[9,1] \rightarrow C[3,1]$, 
both $\kc{Q_1}{K}$ and $ \kc{Q_0}{K}  $ are needed to determine the invariant.

}

\end{remark}

\begin{remark}
{\rm
From the considerations in the preceding remark, it is of interest
whether the number of colorings of a quandle determines that of
another.  There are many such pairs among RIG quandles over all $2977$
knots. 
We note that
many of the pairs do not have epimorphisms between them, and more
studies on this phenomenon may be desirable.
}
\end{remark}





\begin{remark}
{\rm
The quandle $2$-cocycle invariant was defined in \cite{CJKLS} as follows.
Let $\phi: X \times X \rightarrow A$ be a quandle $2$-cocycle, 
for a finite quandle $X$ and an abelian group $A$.
For  a knot  diagram $K$ and  a coloring $x \in \skc{X}{K}$, the weight at a crossing $\tau$ as depicted 
in Figure~\ref{crossing} is defined by $B_\phi (\tau, x) = \epsilon \phi(x, y)$, where $\epsilon = \pm 1$ is the sign of 
$\tau$, defined as $+1$ if the under-arc points right in the figure, otherwise $-1$. 
The $2$-cocycle invariant $\Phi_\phi(K) $ is defined as the multiset 
$\{ \sum_{\tau} B_\phi (\tau, x)  \ | \  x \in \skc{X}{K} \} $. 

The cocycle invariant and the invariant ${\rm Col}_{f} (K) $, which we defined
for an epimorphism $f:  Q_1 \rightarrow Q_0$ for quandles $Q_i$, $i=0,1$,
can be naturally combined as follows.
Let $A$ be an abelian group, and $\phi_i$, $i=0,1$, be $2$-cocycles of $Q_i$.
Then define $\Phi_{f; \phi_1, \phi_0} ( K ) $ to be a multiset 
$$ \{  \   ( \sum_{\tau} B_{\phi_0} (\tau, x) ,   \  \{  \sum_{\tau} B_{\phi_1} (\tau, y )\ | \    y \in f_{\sharp}^{-1} (x) \} 
\ )    \ | \ x \in \skc{Q_0}{K} \} . $$

}
\end{remark}

\section{Applications to other knot invariants}\label{othersec}

In this section we give applications of  the number of quandle
colorings to a few invariants for the 2977 knots in \cite{KI}. 
In particular,
we determine the tunnel number for some of the knots with 11 and 12
crossings in \cite{KI}.   For definitions of knot invariants, we refer
to \cite{KI,Rolf}.


Recall that $\col_X(K)$ denotes the
number of colorings of a knot $K$ by a finite quandle $X$.  Let
$\Lq_X(K)$ denote $ \lceil \log_{|X|} (\col_X(K) ) \rceil $, where $|X|$ denotes the
order  of $X$.

We computed $\Lq_Q(K)$ for the 2977 knots over
a set  ${\cal QL}$ of 439 quandles  consisting of all 431 RIG quandles and 
the 8  quandles $Q_i$, $i=16, 18, 19, 20, 22, 23, 25$ from the 26 quandles $Q_i$ in Table \ref{table-Q}.
The additional 8 quandles are not Alexander.
Let $\mlq(K)$ be the maximum of $\Lq_Q(K)$ over all of the $439$ quandles in ${\cal QL}$.
Among $2977$ knots with $12$ crossings or less, there are 
$1473$ knots $K$ with $\mlq(K)=2$, 
$1441$ knots with $\mlq=3$, and 
$63$ ($15$ 11-crossing and 48 12-crossing) knots with $\mlq=4$.  
Let $\mlq^F(K)$ denote the maximum of $\Lq_X(K)$ over all Alexander RIG quandles of the form 
$\Lambda_p/ ( h(t) )$ where $h(t)$ is irreducible in $\Z_p[t]$ for prime $p$ and different from $t$ and $t-1$. 
%
Note that $\mlq$ and $\mlq^F$ are defined  over all quandles in the  particular set ${\cal QL}$. 
The computational results on these are listed in \cite{BridgeIndex}.



\bigskip

\noindent
{\bf Bridge index}

Let $\bg(K)$ denote the bridge index of a knot $K$.
In \cite{Jozef}, a lower bound of the bridge index is given in terms of  the number of quandle colorings: 
Let $X$ be a finite quandle. 
 For an $n$-bridge presentation of a knot, 
an assignment of colors at the $n$ maxima  determines the colors of the remaining arcs of the diagram, 
if it extends to the entire diagram  and gives a well-defined coloring. 
Hence we have $|\col_X(K)| \leq |X|^n$, where $|X|$ denotes the order of $X$. 
Thus we obtain $\Lq_X(K) \leq \bg(K)$,
and therefore, $\mlq(K)  \leq \bg(K)$.

It is proved in \cite{BZ} 
that a Motesinos knot $K$ with $r$ rational tangle summands has $\bg(K)=r$.
In \cite{Chad}, Montesinos knots were used for determining  
 the bridge indices of knots with $11$ crossings.
 In particular, the set of $4$-bridge 11 crossing knots identified in \cite{Chad} agreed with 
 the list of 11 crossing knots with $\mlq=4$. 
 These are: 
 \begin{eqnarray*}
& &  11a\underline{\ }43, 11a\underline{\ }44, 11a\underline{\ }47, 11a\underline{\ }57, 11a\underline{\ }231, 11a\underline{\ }263, \\
& & 11n\underline{\ }71, 11n\underline{\ }72, 11n\underline{\ }73, 11n\underline{\ }74, 11n\underline{\ }75, 11n\underline{\ }76, 11n\underline{\ }77, 11n\underline{\ }78, 11n\underline{\ }81.
\end{eqnarray*}

At the time of writing, KnotInfo \cite{KI} did not have informalion on the bridge index for 12 crossing knots.
The set of $12$ crossing Montesinos knots with $4$ rational tangle
summands were computed by Slavik Jablan using LinKnot \cite{LinKnot},
and the list was provided to us by Chad Musick. There are 48 knots in
their table.  Our computation showed that there are $48$ knots with 12
crossings with $\mlq=4$, and the list can be found in \cite{BridgeIndex}. 
These two lists coincide.

\bigskip

\noindent
{\bf Nakanishi index}


We use the notations  $\Lambda=\Z[t, t^{-1}]$ and $\Lambda_p=\Z_p[t, t^{-1}]$.
 The Nakanishi index $\NI (K)$ of a knot $K$ 
 is the minimum size of  square presentation matrices of  the Alexander module 
 ($H_1(\tilde{Y})$ as a $\Lambda$-module, where $\tilde{Y}$ is the infinite cyclic covering of the complement $Y$) of $K$ \cite{Naka}.
 The Nakanishi index for $11$ and $12$ crossing knots was blank in \cite{KI} 
 when this project was started (June, 2012).
 We examine how much of the Nakanishi index can be determined for these knots using quandle colorings.

\begin{lemma}
 For any knot $K$ and  Alexander quandle $Q$ of the form 
 $\Lambda_p/ ( h(t) )$ where $h(t)$ is irreducible in $\Z_p[t]$ for prime $p$ and different from $t$ and $t-1$,
 it holds that 
 $$\Lq_Q(K) \leq \NI(K) + 1  \leq \bg(K) . $$
\end{lemma}
\noindent 
{\it Proof.\/} 
As in the proof of Proposition~\ref{dualAlex},
 $\Lambda_p/ (h(t))$ is a finite field $\F_q$ of order $q$ denoted by  $\F(p, h(t))$ in \cite{BF},
  where $q=p^k$ for a prime $p$ and  $k= {\rm deg}(h(t))$. 
  
  It is shown in 
  \cite{BF} that 
if $A$ is a presentation   matrix of the Alexander module,
then the same matrix reduced modulo $p$, $A^{(p)}$, is a presentation matrix of 
  $H_1(\tilde{Y}; \Z_p)$ as a $\Lambda_p$-module.
  Let $\overline{A^{(p)}}$ be the matrix  $A^{(p)}$ with entries reduced modulo $h(t)$. 
  
 Let 
$\psi : \Lambda_p^n \rightarrow \Lambda_p^n $ denote the map corresponding to $A^{(p)}$,
and $\overline{\psi} : \F_q^ n \rightarrow \F_q^n  $ denote the  map on the vector space $\F_q^n$
corresponding to $\overline{A^{(p)}}$. 
Then in \cite{BF} it was shown that the set of colorings $\skc{Q}{K}$ is in one-to-one correspondence 
with the vector space  $ \F_q \oplus {\rm ker}(\overline{\psi} ) $. 

If $\NI (K)=n$, then there is an $n \times n$  presentation matrix $A$ for the Alexander module. 
Then ${\rm ker}(\overline{\psi} )$ has dimension at most $n$. 
Hence $$\kc{Q}{K} = |  \F_q \oplus {\rm ker} (\overline{\psi}) | = q^{1+{\rm dim \,ker}(\overline{\psi})} 
\leq q^{1 + n}, $$
so that we obtain $\Lq_Q (K) \leq 1+n$.
It is well known that $\NI(K) + 1  \leq \bg(K) $. 
$\Box$


 
 \begin{corollary}\label{nakacor}
 For any knot $K$, $\mlq^{F} (K) \leq \NI(K) + 1 \leq \bg(K)$. 
\end{corollary}

For 10 crossing knots, in KnotInfo \cite{KI} as of October 2013, 
the following knots are listed as having $\NI=2$: $10\underline{\ }k$ for 
$$k=61,
63,
65,
69 , 
74,
75,
98,
99,
101, 
103,
122,
123,
140,
142,
144,
155,
157,
160
. $$ 
The list of 10 crossing knots such that Corollary~\ref{nakacor}  determines that $\NI=2$ 
consists of $10\underline{\ }\ell$ for 
$$\ell=61,
63,
65,
74,
75,
98,
99,
103,
115,
122,
123,
140,
142,
144,
155,
157,
163 
. $$ 
Thus Corollary~\ref{nakacor} determines knots with  $\NI=2$ for all but $k=$69, 101, and 160.
We point out that Corollary~\ref{nakacor} implies $\NI=2$
for $k=$115 and 163, these values of NI  are posted in \cite{KI} incorrectly as having $\NI=1$. This error was 
corrected in \cite{Kawauchi} (the knot $10_{163}$ is denoted as $10_{164}$ in \cite{Kawauchi}).

 The  $11$ and $12$ crossing knots with
 $\mlq^{F} (K)=3$, and therefore, $\NI(K)\geq2$ by Corollary \ref{nakacor},
 are listed in \cite{BridgeIndex}. 
 Among these, knots with $\bg(K)=3$ are determined to have $\NI(K)=2$ 
 by Corollary \ref{nakacor}.
  In particular, 11 crossing knots in this list (46 knots) have bridge index 3 by \cite{Chad}, and therefore, 
  have $\NI=2$.

 \bigskip

 \noindent
 {\bf Tunnel number}
 
 Let  $\tau(K)$  denote the tunnel number of a knot $K$.
 At the writing of this article (October 2013), KnotInfo \cite{KI} did not contain 
 the tunnel number for knots with 11 or more crossings.
 
 Let $Q$ be a quandle of the form $\Lambda_p/ (h(t) ) $ for a prime $p$ and an irreducible polynomial $h(t)$ as before.
 It was shown in  \cite{Ishii} (see also \cite{BF}) that 
$\Lq_Q(K)  \leq \tau(K) + 1$ for any knot $K$.
It is well known that $ \tau(K) + 1 \leq \bg(K) $ holds for any knot $K$.
Hence we obtain the following lemma.
\begin{lemma}
If  $\Lq_Q(K)=\bg(K)$ for some quandle $Q$  of the form 
$\Lambda_p/ ( h(t) )$ where $h(t)$ is irreducible in $\Z_p[t]$ for prime $p$ and different from $t$ and $t-1$,
 then 
$\tau(K)=\bg(K)-1$.
\end{lemma}

{}As we mentioned earlier the bridge index was determined for all 11 crossing knots in \cite{Chad}.
For RIG Alexander quandles $Q$  of the form $\Lambda_p/ (h(t) ) $, 
the following 11 crossing knots satisfy the condition $\Lq_Q(K)=\bg(K)=3$,
and therefore they have  tunnel number $2$:
$11a\underline{\ }k $ and $11n\underline{\ } \ell$ where
\begin{eqnarray*}
k &=& 87, 97, 107, 123, 132, 133, 135, 143, 155, 157, 165, 173, 181, 196, 239, 249, \\
& & 277, 291, 293, 297, 314, 317, 321, 322, 329, 332, 340, 347, 352, 354, 366.
\\
\ell &=& 49, 83, 90, 91, 126, 133, 148, 157, 162, 164, 165, 167, 175, 184, 185.
\end{eqnarray*}

These are 46 knots among 552 total of 11 crossing knots. 

 \bigskip

 \noindent
 {\bf Unknotting number}

It is known \cite{Naka} that the unknotting number $u(K)$ of a knot $K$
is bounded below by $\NI(K)$: $\NI(K) \leq u(K)$. We use this fact to
examine the unknotting number of knots in the KnotInfo table.

We consider the knots with  $\mlq^F (K) =3$. 
The knots in the table with $\mlq^F (K) =3$ are posted at 
\cite{BridgeIndex}, 
together with the first simple prime power Alexander quandle  $X$ that gives $\Lq_X(K)=3$.
A file containing 
knots with $\mlq^F (K) =3$ and their unknotting numbers 
is also found at \cite{URL}. 
%
%
 In the list, the notation $[2,3]$ means $u(K)=2$ or $3$,
which we also denote by $u(K)=[2,3]$ below for shorthand.
By Corollary~\ref{nakacor} and  the inequality 
$\NI(K)  \leq u(K)$, for these knots  with  $\mlq^F (K) =3$ 
we obtain $2=\mlq^F (K)-1\leq \NI(K) \leq u(K) $. 
Hence we obtain the following information 
on the Nakanishi index and the unknotting number for these knots with $\mlq^F(K)=3$ and with crossings 11 or 12:

\begin{itemize}

\item[(U1)]
If $u(K)=2$, then it is determined that $\NI(K)=2$. 
The knots that satisfy this condition are listed below.

$11a\underline{ \ }:$ 87,  97, 107, 132, 133, 135, 143, 155, 157, 165, 173, 181, 196, 239, 249, 277, 293, 297, 

\hspace{9.5mm} 314, 317, 321, 322, 332, 347, 352. 

$11n\underline{ \ }:$  90, 164, 175, 185.

$12a\underline{ \ }:$  216, 253, 408, 444, 466, 679, 701, 987, 1183, 1206.

$12n\underline{ \ }:$ 273, 332, 403, 436, 508, 510, 526, 549, 565, 570, 592, 604, 617, 643, 666, 839, 887.

\item[(U2)]
If $u(K)=[1,2]$, then  it is determined that  $\NI(K)=2$ and $u(K)=2$. 
The knots that satisfy this condition are listed below.

$11a\underline{ \ }$: None.

$11n\underline{ \ }$: 49, 83, 91, 157, 162, 165, 167. 

$12a\underline{ \ }$: There are 63 knots in this category. See Appendix~\ref{unk-sec} for the list.

$12n\underline{ \ }$: There are 74 knots in this category. See Appendix~\ref{unk-sec} for the list.

\item[(U3)]
If $u(K)=3$, then we obtain the information that 
$\NI(K)=[2,3]$.  The knots that satisfy this condition are listed below.

$11a\underline{ \ }$: 123. 

$11n\underline{ \ }$: 126, 133, 148, 183.

$12a\underline{ \ }$: 295, 311, 327, 386, 433, 561, 563, 569, 576, 615, 664, 683, 725, 780, 907, 921, 1194.

$12n\underline{ \ }$: 147, 276, 387, 402, 494, 496, 581, 626, 654, 660.

\item[(U4)]
If $u(K)=[2,3]$, then 
$\NI(K)=[2,3]$ and no new information is obtained for $u(K)$.
 The knots that satisfy this condition are listed below.

 $11a\underline{ \ }$: 329.

$11n\underline{ \ }$: None.

$12a\underline{ \ }$: 297, 970, 1286.

$12n\underline{ \ }$: 294, 509, 881.

\item[(U5)]
If $u(K)=[1,2,3]$, then  
$\NI(K)=[2,3]$ and the possibility for the unknotting number is narrowed to $u(K)=[2,3]$. 
The knots that satisfy this condition are listed below.

$11a\underline{ \ }$: None.

$11n\underline{ \ }$: None.

$12a\underline{ \ }$: 244, 291, 376, 381, 481, 493, 494, 634, 886, 1124,  1142, 1202, 1205,  1269, 1288.

$12n\underline{ \ }$: 270, 601, 602, 630, 701, 844, 873.

\item[(U6)]
If $u(K)=[2,3,4]$ ;  $[3,4]$ ; $4$, respectively, 
then $\NI(K)=[2,3,4]$ and  no new information is obtained for $u(K)$.
The knots that satisfy these conditions are listed below.

$11a\underline{ \ }$: 354, 366 ; 291, 340 ; None, respectively.

$11n\underline{ \ }$: None ; None ; None.

$12a\underline{ \ }$: 1097, 1164 ; 973 ; 574, 647.

$12n\underline{ \ }$: None ; 600, 764, 806 ;  386, 518.

\end{itemize}

We observed no other cases in the list.
In particular, in our list of knots with $\mlq^F=3$, there is no knot $K$ with
$u(K)$ containing 5, nor $u(K)=[1,2,3,4]$. All the other positive integer intervals containing less than $5$ 
are listed above.


\begin{remark}
{\rm

Among 12 crossing  knots $K$ with  $\mlq(K)=4$ (recall that there are 48 of them), 
the following 7 knots have $\mlq^F(K)=4$, and 
for each of the 7 knots, the only RIG quandle  $Q$ of type $F(p,h(t))$  such that $Lq_Q(K) = 4$ is
$ Q = C[3,1]$:

$12a\underline{ \ }$: 554, 750.

$12n\underline{ \ }$:  553, 554, 555, 556, and 642.

   These knots have the Nakanishi index at least $3$ from
  Corollary~\ref{nakacor}.
The bridge indices of these knots were determined to be $4$.  
Hence by Corollary~\ref{nakacor}, we obtain that
  $\NI(K)=3$ for these knots.

  The unknotting number for all of these knots except 
  $12n\underline{\ }642$ is listed in \cite{KI} as $[1,2,3]$, and it
  is posted as $[2,3,4]$ for $12n\underline{\ }642$.
 Since $\NI=3$, the unknotting number of these knots except 
  $12n\underline{\ }642$ is determined to be $3$, and it is $[3,4]$
  for $12n\underline{\ }642$.
}
\end{remark}

\noindent
{\it Acknowledgment:}
  We were greatly helped in this project by the construction of all
  connected quandles of order at most 35 by Leandro Vendramin
  \cite{rig}.  We also acknowledge help from Michael Kinyon in showing
  us how to use the program Mace4 \cite{McC} to compute the number of
  coloring of a knot by a quandle. We also thank Larry Dunning for
  showing us how to minimize the number of quandles needed by use of
  binary integer linear programming.  Thanks also to James McCarron
  for his help with the Maple Package Magma and especially with the
  use of the Maple procedure Magma:-AreIsomorphic for determining when
  two quandles are isomorphic.  We would also like to thank David
  Stanovsky for useful information about connected quandles, and Slavik
  Jablan and Chad Musick for valuable information on the bridge index.

  MS was partially supported by  NSF  DMS-0900671 and NIH 1R01GM109459-01. 

\bigskip
\newpage


\appendix

\begin{center}
{\Large\bf Appendix}
\end{center}

\section{Enumeration of quandles and properties of RIG quandles}\label{propertysec}

Sequences enumerating quandles of order $n$ with various properties  may be found in the
Online Encyclopedia of Integer Sequences~\cite{OEIS}. We list  the OEIS identification numbers of some of these sequences below. We also list below some properties and characterizations of RIG quandles that may be found in files at \cite{URL}.

\begin{itemize}

\item
\url{http://oeis.org/A193024}. The number of isomorphism classes of Alexander quandles of order  $n$.

\item
\url{http://oeis.org/A193067}. The number of isomorphism classes of connected Alexander  quandles of order  $n$. 

\item
\url{http://oeis.org/A225744}.
The number of isomorphism classes of connected Generalized Alexander quandles of order  $n$.

\item
\url{http://oeis.org/A196111} The number of isomorphism classes of simple quandles of order  $n$.
A quandle is {\it simple} if it has more than one element, and if it has no homomorphic images other than itself or the singleton quandle.  However, in \cite{AG}, it is further assumed that simple quandles are not trivial. This also rules out the quandles of order 2. 

\item
\url{http://oeis.org/A177886} The number of isomorphism classes of latin quandles of order  $n$.
A {\it Latin quandle} is a quandle such that for each $a \in X$, the left translation ${\cal L}_a$ is a bijection. That is, the multiplication table of the quandle is a Latin square.

\item
\url{https://oeis.org/A226172} The number of isomorphism classes of faithful connected quandles of order  $n$. A quandle is {\it faithful} if the mapping $a \mapsto {\cal R}_a$ from $X$ to ${\rm Inn}(X)$ is an injection.

\item
\url{https://oeis.org/A226173} The number of isomorphism classes of connected keis (involutory quandles) of order  $n$.
A quandle $X$ is {\it involutory}, or a {\it kei}, 
if the right translations  are  involutions: ${\cal R}_a^2 ={\rm  id}$ for all $a \in X$. 

\item
\url{https://oeis.org/A226174} The number of isomorphism classes of self-dual connected quandles of order $n$.

\end{itemize}

\noindent We also obtained the following information about RIG quandles. This is available at \cite{URL}: 

\begin{itemize}
\item
Subquandles of RIG quandles that are RIG quandles and trivial quandles.
We also identified the orders of subquandles that are neither connected nor trivial.

\item 
The connected RIG Alexander quandles are presented  as 
$\Z[t] /J$ where $J =  (  f(t) ) $  or $J = (f(t),g(t))$ except for one Alexander RIG quandle, namely $C[27,17]$,  which is not a cyclic $\Z[t,t^{-1}]$-module. 

\item
We identified RIG quandles that are  products of other RIG quandles.

\item
We identified dual quandles of RIG quandles. In particular, if the dual is itself, 
then it is self-dual. 

\item
A function $\phi: X \times X \rightarrow A$ for an abelian group $A$ is 
called a {\it quandle $2$-cocycle} if it satisfies 
$$ \phi (x, y)- \phi(x,z)+ \phi(x*y, z) - \phi(x*z, y*z)=0$$
and $\phi(x,x)=0$ for any $x,y,z \in X$.
For a quandle $2$-cocycle $\phi$, 
$X \times A$ becomes a quandle 
by $(x, a) * (y, b)=(x*y, a+\phi(x,y))$ for $x, y \in X$, $a,b \in A$,
and it is called an {\it abelian extension} of $X$ by $A$.
See \cite{CENS}, for example, for more details of extensions.
A list of RIG quandles that are abelian extensions of RIG quandles is found in 
\cite{URL}. 

\item
For each RIG quandle $C[n,i]$ we identified
a set Gen$[n,i] $ which is a smallest set of generators for $C[n,i]$.
We also produced the set of all $[n,i]$ such that any two elements of $C[n,i]$ will generate $C[n,i]$. 
Note that if $[n,i] $ is in this set  then there are no non-trivial subquandles of $C[n,i]$.

\end{itemize}

\begin{remark}
{\rm
It is known that connected quandles of prime order $p$ and those of order $p^2$ 
are Alexander quandles. 
By examining characterizations and properties of RIG quandles, 
listed above, one might wonder if these results might generalize. 
For example, RIG quandles of order $2^n$ for some $n$, and those with order $p q$ for primes 
$p, q > 3$ are Alexander. 
However, inspections of higher order conjugation quandles reveal that
there 
 are quandles with order $64$ and $55$  that are not Alexander. 
 The former is  a size 64 conjugacy 
class of GAP's SmallGroup$(1344, 11699)$.  
Two of the latter 
are
a conjugacy class of PSL$(2,11)$  and  a conjugacy class  of  SmallGroup$( 1210, 7 )$.

}
\end{remark}

\begin{remark}
{\rm
As in \cite{HMN} and \cite{rig}, for computational purposes we
represent quandles by matrices. The entries of the matrix $A$ of a
quandle of order $n$ are integers $1,\ldots,n$ and the quandle
operation is given by $i * j = A_{i,j}$.  We note that the quandles in
the RIG package are left distributive whereas for knot colorings we
prefer right distributive quandles. Hence we have reversed the order
of the product in the RIG quandle by taking the transposes of the matrices
representing the RIG quandles.  In the RIG package the $i$-th quandle
of order $n$ is denoted by SmallQuandle$(n,i)$, for which we use the
notation $C[n,i]$.

}
\end{remark}


\clearpage


\section{Tables of  quandles used to distinguish knots and their mirrors}\label{tablesec}


\begin{appentable}{\rm List of $26$ quandles 
 that distinguish $2977$ knots up to orientation and mirror image.  The quandle matrices can be  found at   \cite{URL}.}
\label{table-Q}

\begin{center}
\begin{tabular}{|c|c|ll|}
\hline
Quandle  & order         & GAP ID for Inn($Q$) & StructureDescription(Inn($Q$))\\
\hline
$Q_{1}$ & $12$ & SmallGroup$(60,5)$ & $A_5$\\ 
$Q_{2}$ & $13$ & SmallGroup$(52,3)$ & $C_{13} : C_4$\\ 
$Q_{3}$ & $13$ & SmallGroup$(52,3)$ & $C_{13} : C_4$\\ 
$Q_{4}$ & $13$ & SmallGroup$(156,7)$ & $(C_{13} : C_4) : C_3$\\ 
$Q_{5}$ & $13$ & SmallGroup$(156,7)$ & $(C_{13} : C_4) : C_3$\\ 
$Q_{6}$ & $15$ & SmallGroup$(60,5)$ & $A_5$\\ 
$Q_{7}$ & $17$ & SmallGroup$(136,12)$ & $C_{17} : C_8$\\ 
$Q_{8}$ & $17$ & SmallGroup$(136,12)$ & $C_{17} : C_8$\\ 
$Q_{9}$ & $20$ & SmallGroup$(120,34)$ & $S_5$\\ 
$Q_{10}$ & $24$ & SmallGroup$(168,42)$ & $PSL(3,2)$\\ 
$Q_{11}$ & $25$ & SmallGroup$(75,2)$ & $(C_5 \times C_5) : C_3$\\ 
$Q_{12}$ & $27$ & SmallGroup$(702,47)$ & $((C_3 \times C_3 \times C_3) : C_{13}) : C_2$\\ 
$Q_{13}$ & $27$ & SmallGroup$(702,47)$ & $((C_3 \times C_3 \times C_3) : C_{13}) : C_2$\\ 
$Q_{14}$ & $30$ & SmallGroup$(120,34)$ & $S_5$\\ 
$Q_{15}$ & $32$ & SmallGroup$(160,199)$ & $((C_2 \times Q_8) : C_2) : C_5$\\ 
$Q_{16}$ & $40$ & SmallGroup$(320,1635)$ & $((C_2 \times C_2 \times C_2 \times C_2) : C_5) : C_4$\\ 
$Q_{17}$ & $40$ & SmallGroup$(320,1635)$ & $((C_2 \times C_2 \times C_2 \times C_2) : C_5) : C_4$\\ 
$Q_{18}$ & $42$ & SmallGroup$(168,42)$ & $PSL(3,2)$\\ 
$Q_{19}$ & $48$ & SmallGroup$(288,1025)$ & $(A_4 \times A_4) : C_2$\\ 
$Q_{20}$ & $60$ & SmallGroup$(660,13)$ & $PSL(2,11)$\\ 
$Q_{21}$ & $72$ & SmallGroup$(576,8652)$ & $(A_4 \times A_4) : C_4$\\ 
$Q_{22}$ & $72$ & SmallGroup$(504,156)$ & $PSL(2,8)$\\ 
$Q_{23}$ & $84$ & SmallGroup$(1512,779)$ & $PSL(2,8) : C_3$\\ 
$Q_{24}$ & $84$ & SmallGroup$(1512,779)$ & $PSL(2,8) : C_3$\\ 
$Q_{25}$ & $90$ & SmallGroup$(720,765)$ & $A_6 . C_2$\\ 
$Q_{26}$ & $182$ & SmallGroup$(1092,25)$ & $PSL(2,13)$\\ \hline
\end{tabular}
\end{center}

\end{appentable}





\clearpage

\begin{appentable}{\rm List of $23$ quandles
 that distinguish $1058$ knots from their mirror images. 
 The quandle matrices can be  found at   \cite{URL}.}
\label{table-T}

\begin{center}
\begin{tabular}{|c|c|ll|}
\hline
Quandle  & order         & GAP ID for Inn($Q$) & StructureDescription(Inn($Q$))\\
\hline
$T_{1}$ & $351$ & SmallGroup$(2106,0)$ & $(((C_3 \times C_3 \times C_3) : C_{13}) : C_3) : C_2$\\ 
$T_{2}$ & $504$ & SmallGroup$(3024,0)$ & $C_2 \times (PSL(2,8) : C_3)$\\ 
$T_{3}$ & $504$ & SmallGroup$(4536,0)$ & $PSL(2,8) : C_9$\\ 
$T_{4}$ & $18$ & SmallGroup$(216,90)$ & $(((C_2 \times C_2) : C_9) : C_3) : C_2$\\ 
$T_{5}$ & $27$ & SmallGroup$(216,86)$ & $((C_3 \times C_3) : C_3) : C_8$\\ 
$T_{6}$ & $27$ & SmallGroup$(486,41)$ & $((C_3 . ((C_3 \times C_3) : C_3) = (C_3 \times C_3) . (C_3 \times C_3)) : C_3) : C_2$\\ 
$T_{7}$ & $28$ & SmallGroup$(168,43)$ & $((C_2 \times C_2 \times C_2) : C_7) : C_3$\\ 
$T_{8}$ & $720$ & SmallGroup$(7920,0)$ & $M11$\\ 
$T_{9}$ & $112$ & SmallGroup$(1344,816)$ & $(((C_2 \times C_2 \times C_2) . (C_2 \times C_2 \times C_2)) : C_7) : C_3$\\ 
$T_{10}$ & $112$ & SmallGroup$(1344,816)$ & $(((C_2 \times C_2 \times C_2) . (C_2 \times C_2 \times C_2)) : C_7) : C_3$\\ 
$T_{11}$ & $117$ & SmallGroup$(1053,51)$ & $((C_3 \times C_3 \times C_3) : C_13) : C_3$\\ 
$T_{12}$ & $162$ & SmallGroup$(1296,2890)$ & $(C_3 . (((C_3 \times C_3) : Q_8) : C_3) = (((C_3 \times C_3) : C_3) : Q_8) . C_3) : C_2$\\ 
$T_{13}$ & $162$ & SmallGroup$(1296,2891)$ & $((((C_3 \times C_3) : C_3) : Q_8) : C_3) : C_2$\\ 
$T_{14}$ & $192$ & SmallGroup$(1344,814)$ & $(C_2 \times C_2 \times C_2) . PSL(3,2)$\\ 
$T_{15}$ & $125$ & SmallGroup$(1500,36)$ & $(((C_5 \times C_5) : C_5) : C_4) : C_3$\\ 
$T_{16}$ & $135$ & SmallGroup$(1620,421)$ & $((C_3 \times C_3 \times C_3 \times C_3) : C_5) : C_4$\\ 
$T_{17}$ & $135$ & SmallGroup$(1620,421)$ & $((C_3 \times C_3 \times C_3 \times C_3) : C_5) : C_4$\\ 
$T_{18}$ & $168$ & SmallGroup$(1512,779)$ & $PSL(2,8) : C_3$\\ 
$T_{19}$ & $64$ & SmallGroup$(448,179)$ & $((C_2 \times C_2 \times C_2) . (C_2 \times C_2 \times C_2)) : C_7$\\ 
$T_{20}$ & $64$ & SmallGroup$(768,1083508)$ & $((((C_8 \times C_4) : C_2) : C_2) : C_2) : C_3$\\ 
$T_{21}$ & $64$ & SmallGroup$(768,1083509)$ & $((((C_8 \times C_4) : C_2) : C_2) : C_2) : C_3$\\ 
$T_{22}$ & $40$ & SmallGroup$(320,1635)$ & $((C_2 \times C_2 \times C_2 \times C_2) : C_5) : C_4$\\ 
$T_{23}$ & $84$ & SmallGroup$(1512,779)$ & $PSL(2,8) : C_3$\\ 
\hline
\end{tabular}
\end{center}

\end{appentable}



\clearpage

\section{Knots with unknotting number 2}\label{unk-sec}

In this section we list knots with 12 crossings such that their unknotting number is listed as 
1 or 2 in KnotInfo~\cite{KI}, and our computation that $\mlq^F=3$ determines that their unknotting number is 
in fact 2. 

\bigskip

The list of alternating knots is $12a\underline{ \ }k $ for $k=$

\noindent
 100,
 177,
 215,
 218,
 245,
 248,
 249,
 265,
 270,
 279,
 298,
 312,
 332,
 347,
 348,
 396,
 413,
 427,
 429,
 435,
 448,
 465,
 475,
 503,
 594,
 703,
 712,
 742,
 769,
 787,
 806,
 808,
 810,
 868,
 873,
 895,
 904,
 905,
 906,
 941,
 949,
 960,
 975,
 990,
 1019,
 1022,
 1026,
 1053,
 1079,
 1092,
 1093,
 1102,
 1105,
 1123,
 1152,
 1167,
 1181,
 1225,
 1229,
 1251,
 1260,
 1280,
 1283.

\bigskip

The list of non-alternating knots is $12n\underline{ \ }\ell $ for $\ell=$

\noindent
144,
145,
257,
268,
269,
274,
297,
333,
334,
355,
356,
357,
379,
380,
388,
389,
393,
394,
397,
414,
420,
440,
442,
460,
462,
480,
495,
498,
505,
533,
546,
567,
571,
582,
583,
598,
605,
611,
622,
636,
637,
651,
652,
669,
706,
714,
717,
737,
742,
745,
746,
752,
756,
757,
760,
779, 
781,
798,
813,
817,
837,
838,
840,
843,
846,
847,
869,
874,
876,
877,
878,
879,
883.

\newpage


\section{Remarks on the computation of $\kc{Q}{K} $} \label{computationsec}

We used several techniques to compute $\kc{Q}{K}$ for quandles $Q$ and
knots $K$. We take for each knot the braid representation given at
KnotInfo \cite{KI}.  We take the braids to be oriented from top to
bottom. This induces an orientation of the knots.  The two parameters
that are most important for computing the number of colorings are the
braid index, $b(K)$ (the number of strands) and the order $|Q|$ of the
quandle. A straightforward computation requires $|Q|^{b(K)}$ steps.
Using the fact that our quandles are connected, i.e., the group
$Inn(Q)$ acts transitively on $Q$ we may assume that one of the arcs
has a fixed color. This allowed us to reduce the number of steps to
$|Q|^{b(K)-1}$, which reduces the time by a factor of $1/|Q|$. In the
case of quandles of large order this makes a noticeable difference in
program running times (hours vs months).

We initially used Mace4 \cite{McC} to find the number of
colorings. Later on, we wrote a program in C to verify that the
results we obtained using Mace4 were correct. 
The C
program uses the same 
braid representations as those used for
the Mace4 computations.

Let $K$ be a knot with a braid index $b(K)$.
Let $Q$ be a quandle of order $n$. Let $a_i$,
$i=1, \ldots, b$, represent the element of a given quandle assigned to the top of the $i$-th 
strand of a braid, where $b=b(K)$ is the braid index.
 To count the number of knot
colorings, we
find the colorings at the bottom strands of the braid,
and check assignments such that the top and bottom colors are the same for each strand,
compare Figure~\ref{figure8}.
The following pseudo code demonstrates that our algorithm runs in
$O(n^{b-1})$ time, $n=|Q|$. \\

\lstset{language=C}
\begin{lstlisting}[frame=single]
/* fix strand a1 to be 0. This allows us to get away *
 * with coloring b-1 strands instead of b strands.   */
 a1=0;
 for(a2=0; a2 < n; a2++) {
   for(a3 ...) {
     for(a4 ...) {
       for(a5 ...) {
         ... up to the braid index 'b'
           /* constant time quandle coloring of knot strands */
       }
     }
   }
 }
\end{lstlisting}


We implemented a multithreaded version of the program using POSIX
threads 
to further improve performance. Combined with
fixing a color of one strand, 
we found that multithreading reduces the
computation time considerably on multicore processors. 
The unexpectedly large speedup of
multithreading over trivial parallelization may be due to operating
system characteristics and prevention of CPU cache swapping to the
slower RAM when the operating system is switching tasks.

The computations were performed on a GNU/Linux computer equipped with
24 Gigabytes of RAM and an Intel Core i7 processor with 4 cores. The
kernel was compiled with a task switching latency of 100hz.

\begin{figure}[H]
  \begin{center}
    \includegraphics[height=4in]{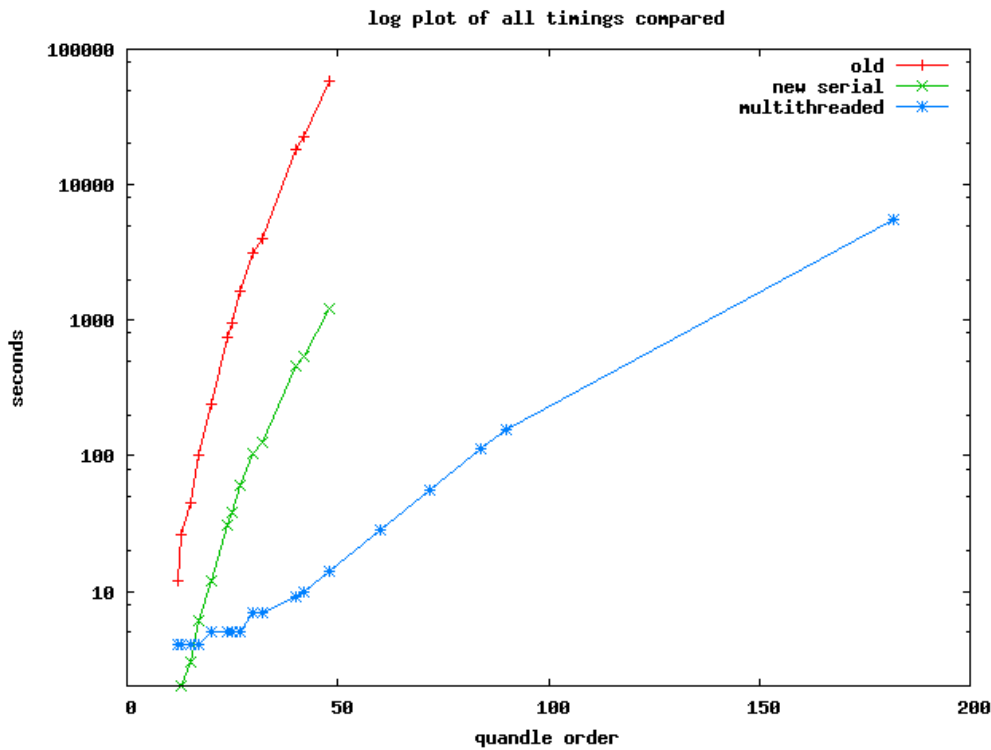}
    \caption{Log plot of timings}\label{graph}
  \end{center}
\end{figure}

In Table~\ref{table-comp}, the first column indicates $26$ quandles in Table \ref{table-Q},
the second column gives  the order of the  quandle, the third column is the computation time 
for all 2977 knots with at most 12 crossings
using  a single thread, the fourth column is time using multiple threads, and the  last column is 
the ratio of column $3$ and column $4$.
A log plot of timings of three different coloring programs is given in  Figure~\ref{graph}.

\begin{table}
\begin{center}
\begin{tabular}{|c|c|r|l|r|}
\hline
quandle    & order &  serial seconds &  multithread seconds & speedup factor\\
\hline
\hline
$Q_{1}$    & 12    &  1              &  4                   & 0.25 $\times$ (slower)\\
$Q_{2}$    & 13    &  2              &  4                   & 0.5 $\times$ (slower)\\
$Q_{3}$    & 13    &  2              &  4                   & 0.5 $\times$ (slower)\\
$Q_{4}$    & 13    &  2              &  4                   & 0.5 $\times$ (slower)\\
$Q_{5}$    & 13    &  2              &  4                   & 0.5 $\times$ (slower)\\
$Q_{6}$    & 15    &  3              &  4                   & 0.75 $\times$ (slower)\\
\hline
$Q_{7}$    & 17    &  6              &  4                   & 1.5 $\times$ (faster)\\
$Q_{8}$    & 17    &  6              &  4                   & 1.5 $\times$ (faster)\\
$Q_{9}$    & 20    &  12             &  5                   & 2.4 $\times$ (faster)\\
$Q_{10}$   & 24    &  31             &  5                   & 6.2 $\times$ (faster)\\
$Q_{11}$   & 25    &  38             &  5                   & 7.6 $\times$ (faster)\\
$Q_{12}$   & 27    &  60             &  5                   & 12 $\times$ (faster)\\
$Q_{13}$   & 27    &  60             &  5                   & 12 $\times$ (faster)\\
$Q_{14}$   & 30    &  105            &  7                   & 15 $\times$ (faster)\\
$Q_{15}$   & 32    &  125            &  7                   & 17.86 $\times$ (faster)\\
$Q_{16}$   & 40    &  455            &  9                   & 50.56 $\times$ (faster)\\
$Q_{17}$   & 40    &  455            &  9                   & 50.56 $\times$ (faster)\\
$Q_{18}$   & 42    &  542            &  10                  & 54.2 $\times $ (faster)\\
$Q_{19}$   & 48    &  1215           &  14                  & 86.79 $\times $ (faster)\\
$Q_{20}$   & 60    &  (hours)        &  28                  & ? \\
$Q_{21}$   & 72    &  (days)         &  56                  & ? \\
$Q_{22}$   & 72    &  (days)         &  56                  & ? \\
$Q_{23}$   & 84    &  (days)         &  112                 & ? \\
$Q_{24}$   & 84    &  (days)         &  112                 & ? \\
$Q_{25}$   & 90    &  (weeks)        &  155                 & ? \\
$Q_{26}$   & 182   &  (months)       &  5602                & (a lot faster) \\
\hline
\end{tabular}
\caption{Computational time comparison}\label{table-comp}
\end{center}
\end{table}


\bigskip

Multithreading was used to perform computations in parallel for the
outermost loop of our knot coloring program. We set the number of
threads to be the order of the quandle. None of the child threads
communicate with each other. The parent thread waits for all of the
child threads to return the number of colorings from its search. The
parent sums each of the results from the child threads. The number of
colorings of the knot is the sum of the counts from the child threads.
The multithreaded version of the program is required to compute the
colorings by  the quandle of order 182.

\bigskip

\begin{lstlisting}[frame=single]
/* fix strand a1 to be 0. This allows us to get away *
 * with coloring b-1 strands instead of b strands.   */
 a1=0;

 /* when multithreading, we compute b-2 strands since a2 is fixed */
 a2=thread_id;
 /* for each thread */
 for(a3=0; a3 < n; a3++) {
   for(a4 ...) {
     for(a5 ...) {
       ... up to the braid index 'b'
         /* constant time quandle coloring of knot strands */
     }
   }
 }
\end{lstlisting}

The output of our program is the number ${\rm Col}_{Q}(K)_0$ of 
colorings with the color of the first strand fixed. 
The result is then converted to the 
number of non-trivial colorings ${\rm Col}^{N}_{Q}(K)$ using 
the formula ${\rm Col}^{N}_{Q}(K)
= |Q|({\rm Col}_{Q}(K)_0-1)$. 


All of the programs we used are provided on-line at \cite{URL
} 
with source and data sets
available. We made numerous checks, but would be pleased if our
computations can be independently confirmed.


\end{document}